\documentclass[11pt,a4paper,draft]{article}
\usepackage{mathrsfs}

\usepackage[leqno]{amsmath}
\usepackage{amssymb,amsthm,upref,amscd}
\usepackage[T1]{fontenc}
\numberwithin{equation}{section}
\input scrload.tex

\theoremstyle{plain}
\newtheorem{Thm}{Theorem}[section]
\newtheorem{Lem}[Thm]{Lemma}
\newtheorem{Cor}[Thm]{Corollary}
\newtheorem{Prop}[Thm]{Proposition}
\newtheorem*{Thm*}{Theorem}
\theoremstyle{definition}

\newtheorem{Rem}[Thm]{Remark}

\newcommand{\equ}{equation}
\newcommand{\C}{\mathbb{C}}

\newcommand{\R}{\mathbb{R}}

 \DeclareMathOperator{\dist}{dist}

 \DeclareMathOperator{\sgn}{sgn}

\let\nhatoksa=\theenumi
\let\nhatoksb=\labelenumi
\let\nhatoksc=\theenumii
\let\nhatoksd=\labelenumii
\newlength{\nhalengtha}
\newlength{\nhalengthb}
\newlength{\nhalengthc}
\newcommand{\resetenum}{
\let\theenumi=\nhatoksa
\let\labelenumi=\nhatoksb
\let\theenumii=\nhatoksc
\let\labelenumii=\nhatoksd
\setlength{\leftmargini}{\nhalengtha}
\setlength{\leftmarginii}{\nhalengthb}
\setlength{\labelwidth}{\nhalengthc}
}

\newcommand\cd{\mathcal{D}}

\newcommand\ci{\mathcal{I}}
\newcommand\cj{\mathcal{J}}

\newcommand\cm{\mathcal{M}}
\newcommand\cn{\mathcal{N}}

\newcommand\ct{\mathcal{T}}

\def\msc{\mathscr{C}}

\def\msf{\mathscr{F}}
\def\msg{\mathscr{G}}

\def\msj{\mathscr{J}}
\def\msk{\mathscr{K}}
\def\msl{\mathscr{L}}
\def\msm{\mathscr{M}}
\def\msn{\mathscr{N}}

\def\mss{\mathscr{S}}
\def\mst{\mathscr{T}}

\def\pa {\partial}

\def\op{\oplus}

\def\De{\Delta}

\def\ka {\kappa}
\def\al {\alpha}
\def\bt {\beta}
\def\de {\delta}
\def\Ga {\Gamma}
\def\ga {\gamma}
\def\lm {\lambda}
\def\Lam {\Lambda}
\def\om {\omega}

\def\sa {\sigma}
\def\vr {\varepsilon}
\def\va {\varphi}

\def\sgn{\hbox{sgn}}

\newcommand{\bkt}[1]{\left(#1\right)}

\newcommand{\norm}[1]{\left\|#1\right\|}
\newcommand{\inp}[2]{\left\langle#1,#2\right\rangle}
\newcommand{\jdz}[1]{\left|#1\right|}

\title{On the concentration of semi-classical states for a nonlinear
Dirac-Klein-Gordon system}

\author{Yanheng Ding, \, Tian Xu   \and {\small Institute
of Mathematics, AMSS, Chinese Academy of Sciences,} \\
{\small  100190 Beijing, China} }

\date{}

\begin{document}
\maketitle

\begin{abstract}
In the present paper, we study the semi-classical approximation of
a Yukawa-coupled massive Dirac-Klein-Gordon system
with some general nonlinear self-coupling. We prove
that for a constrained coupling constant there exists
a family of ground states of the semi-classical problem,
for all $\hbar$ small, and show that the family concentrates
around the maxima of the nonlinear potential as $\hbar\to0$.
Our method is variational and relies upon a delicate
cutting off technique. It allows us to overcome the
lack of convexity of the nonlinearities.

\vspace{.5cm}


\vspace{.5cm}


\end{abstract}

\section{Introduction and main result}

In this paper we study the solitary wave solutions of
the massive Dirac-Klein-Gordon system involving an external
self-coupling:
\begin{\equ}\label{M-D0}
\left\{
\begin{aligned}
&i\frac\hbar c\,\pa_t \psi
+ i\hbar\sum_{k=1}^3\al_k\pa_k \psi-mc\bt\psi-\lm\phi\bt\psi=f(x,\psi)\\
&\frac{\,\hbar^2}{\,c^2}\,\pa_t^2\phi
- \hbar^2 \De \phi +M \phi=4\pi\lm(\bt \psi)\cdot\psi
\end{aligned}  \right.
\end{\equ}
for $(t,x)\in\R\times\R^3$, where $c$ is the speed of light, $\hbar$
is Planck's constant, $\lm>0$ is coupling constant, $m$ is
the mass of the electron and $M$ is the mass of the meson (we use
the notation $u \cdot v$ to express the inner product of $u,
v\in\C^4$). Here
$\alpha_1$, $\alpha_2$, $\alpha_3$ and $\beta$ are $4\times4$
complex Pauli matrices:
\[
\beta=\left(
\begin{array}{cc}
I&0\\
0&-I
\end{array}\right), \quad \alpha_k=\left(
\begin{array}{cc}
0&\sigma_k\\
\sigma_k&0
\end{array}\right), \quad k=1,2,3,
\]
with
\[
\sigma_1=\left(
\begin{array}{cc}
0&1\\
1&0
\end{array}\right), \quad \sigma_2=\left(
\begin{array}{cc}
0&-i\\
i&0
\end{array}\right), \quad \sigma_3=\left(
\begin{array}{cc}
1&0\\
0&-1
\end{array}\right)  .
\]
System \eqref{M-D0} arises in mathematical models of particle physics,
especially in nonlinear topics. Physically, system \eqref{M-D0}
describes the Dirac and Klein-Gordon equations coupled through the
Yukawa interaction between a Dirac field $\psi\in\C^4$
and a scalar field $\phi\in\R$ (see \cite{Bjorken}).
This system is inspired by approximate
descriptions of the external force involve only functions of fields.
The nonlinear self-coupling $f(x,\psi)$, which describes
a self-interaction in Quantum electrodynamics,
gives a closer description of many particles found
in the real world. Various
nonlinearities are considered to be
possible basis models for unified field theories (see \cite{FLR},
\cite{FFK}, \cite{Iva} etc. and references therein).

System \eqref{M-D0} with null external self-coupling, i.e., $f\equiv 0$, has
been studied for a long time and results are available concerning
the Cauchy problem (see \cite{Bournaveas}, \cite{Chadam1},
\cite{Chadam3}, \cite{Machihara}, \cite{Selberg}
 etc.).
The first result on the global existence and
uniqueness of solutions of \eqref{M-D0} (in one space dimension)
was obtained by J. M. Chadam
in \cite{Chadam1} under suitable assumptions on the initial data.
For later developments, we mention, e.g.,
that J. M. Chadam and Robert T. Glassey \cite{Chadam3} yield the
existence of a global solution in three space dimensions.
In \cite{Bournaveas}, N. Bournaveas obtained
low regularity solutions of the Dirac-Klein-Gordon system
by using classical Strichartz-type time-space estimates.

As far as the existence of stationary solutions (solitary
wave solutions) of \eqref{M-D0} is concerned,
there is a pioneering work by
M. J. Esteban, V. Georgiev and E. S\'{e}r\'{e} (see \cite{Sere1})
in which a multiplicity result is studied.
Here, by stationary solution, we mean a solution of the type
\begin{\equ}\label{stationary}
\left\{
\aligned
&\psi(t,x)=\va(x)e^{-i\xi t/\hbar}, \quad \xi\in\R, \quad \va:\R^3\to \C^4, \\
&\phi=\phi(x) \, .
\endaligned \right.
\end{\equ}
In \cite{Sere1}, using the variational arguments,
the authors obtained infinitely many solutions for $\xi\in(-\frac{mc}{\hbar},0)$
under the assumption
\[
\va(x)=\left(
\begin{array}{cc}
v(r) &  \Big( \begin{array}{c} 1 \\ 0 \end{array} \Big) \\
iu(r) & \Big( \begin{array}{c} \cos\vartheta \\ e^{i\tau}\sin\vartheta
               \end{array} \Big) \\
\end{array}
\right)
\]
where $(r,\vartheta,\tau)$ are the spherical coordinates of $x\in\R^3$.

We emphasize that the works mentioned above were mainly concerned  with
the autonomous system with null self-coupling. Besides, limited
work has been done in the semi-classical approximation.
In the present paper we are devoted to the existence and
concentration phenomenon of stationary
semi-classical solutions to system \eqref{M-D0}.
For small $\hbar$, the solitary waves are
referred to as semi-classical states. To describe the transition
from quantum to classical mechanics, the existence of solutions
$(\va_\hbar,\phi_\hbar)$, $\hbar$ small,
possesses an important physical interest.
More precisely, for ease of notations, denoted by $\vr=\hbar$,
$\al=(\al_1, \al_2, \al_3)$ and
$\al\cdot\nabla=\sum_{k=1}^3\al_k\pa_k$, we are concerned
with (substitute \eqref{stationary} in \eqref{M-D0})
the following stationary nonlinear Dirac-Klein-Gordon system:
\begin{\equ}\label{D2}
\left\{
\begin{aligned}
&i\vr\,\al\cdot\nabla \va-a\bt \va + \om \va
- \lm \phi\bt\va = W(x)g(|\va|) \va  \, ,\\
&-\vr^2\De \phi +M \phi =4\pi \lm (\bt\va)\cdot\va  \, .
\end{aligned}\right.
\end{\equ}
where $a=mc>0$ and $\om\in\R$.

On the nonlinear self-coupling, writing $G(|w|):=\int^{|w|}_0g(s)sds$,
we make the following hypotheses:
\begin{itemize}
\item[$(P_0)$] $W\in C(\R^3)\cap L^\infty(\R^3)$
 with $\inf W>0$ and $\limsup\limits_{|x|\to\infty}W(x)<\max W(x)$;

\item[$(G_1)$] {\it $g(0)=0$, $g\in C^1(0,\infty)$, $g'(s)>0$ for $s>0$, and
there exist $p\in (2, 3)$, $c_1>0$ such that $g(s)\leq
c_1(1+s^{p-2})$ for $s\geq 0$}\ ;

\item[$(G_2)$]{\it  there exist $\sigma>2$, $\theta>2$ and $c_0>0$ such that
$c_0s^\sigma\leq G(s)\leq \frac{1}{\theta}g(s)s^2$ for all $s>0$}\ .
\end{itemize}
A typical example is the power function $g(s)=s^{\sigma-2}$.

For showing the concentration phenomena,  we set $m:=\max_{x\in\R^3}W(x)$ and
$$
\msc:=\{x\in\R^3: \, W(x)=m\}.
$$
Our result reads as

\begin{Thm}\label{main theorem}
Assume that $\om\in(-a,a)$, $(P_0)$ and $(G_1)$-$(G_2)$  are
satisfied. Then there exists $\lm_0>0$ such that
given $\lm\in(0,\lm_0]$, for all $\vr >0$ small,
\begin{enumerate}
\item      The system \eqref{D2} possesses at least one ground state solution
$(\va_\vr,\phi_\vr)\in \cap_{q\geq2}W^{1,q}(\R^3,\C^4) \times C^2(\R^3,\R)$.

\item      The set of ground state solutions is compact in
$H^1(\R^3,\C^4) \times H^1(\R^3,\R)$.

\item      If additionally $\nabla W$ is bounded, then
\begin{itemize}
\item[$(i)$] There is a maximum point $x_\vr $ of
$|\va_\vr|$ with
$\lim_{\vr \to0}\dist (x_\vr ,\msc)=0$ such that the pair $(u_\vr, V_\vr)$,
where $u_\vr (x):=\va_\vr (\vr  x+x_\vr )$ and $V_\vr:=\phi_\vr (\vr  x+x_\vr )$,
converges in $H^1\times H^1$ to a ground state solution of (the  limit
equation)
\begin{\equ}\label{the limit problem}
\left\{
\aligned
&i \al\cdot\nabla u-a\bt u + \om u
- \lm V \bt u = mg(|u|) u  \\
&- \De V +M \, V =4\pi \lm (\bt u)\cdot u
\endaligned \right. \, .
\end{\equ}
\item[$(ii)$]
$ \jdz{\va_\vr (x)}\leq C
\exp{\big(-\frac{\,c}{\,\vr } |x-x_\vr| \big)} $
for some $C,c>0$.
\end{itemize}
\end{enumerate}
\end{Thm}

It is standard that (\ref{D2}) is equivalent to, by letting $u(x)=
\va(\vr x)$ and $V(x)=\phi(\vr x)$,
\begin{\equ}\label{D3}
\left\{
\begin{aligned}
&i \al\cdot\nabla u - a\bt u + \om u - \lm V\bt u = W_\vr(x) g(|u|)u\\
&-\De V + M \, V = 4\pi \lm  (\bt u)\cdot u
\end{aligned}\right.
\end{\equ}
where $W_\vr (x)=W(\vr x)$. We will in the sequel focus on this equivalent
problem. Our proofs are variational: the
semiclassical solutions that are obtained as critical
points of an energy functional $\Phi_\vr $ associated to the
equivalent problem \eqref{D3}.

There have been a large number of works on existence
and concentration phenomenon of semi-classical
states of nonlinear Schr\"odinger-Poisson
systems arising in the \textit{non-relativistic} quantum mechanics,
see, for example, \cite{Ambrosetti1, Ambrosetti2, Azzollini}
and their references. And, only very recently, the papers
\cite{Tian, Tian2}
studied the existence of a family of semi-classical
ground states of Maxwell-Dirac system and showed
that the family concentrates around some certain sets as $\vr\to0$.
It is quite natural to ask if certain similar
results can be obtained for nonlinear Dirac-Klein-Gordon
systems arising in the relativistic quantum mechanics.
Mathematically, the problems in Dirac-Klein-Gordon systems are difficult because
they are strongly indefinite in the sense that both the negative and positive
parts of the spectrum of Dirac operator are unbounded
and consist of essential spectrums.
%As far as the authors known there have been no results
%on the existence and concentration phenomenon of semiclassical solutions to
%nonlinear Dirac-Klein-Gordon systems.

It should be pointed out that Ding, jointly with co-authors,
developed some technique arguments to obtain the existence
and concentration of semi-classical solutions for
nonlinear Dirac equations (not for Dirac-Klein-Gordon
system), see \cite{Ding2010, Ding2012, Ding Ruf}.
Compared with
the papers, difficulty arises in the Dirac-Klein-Gordon system because of
the presence of the action for a meson field $\phi$. In order
to overcome this obstacle, we develop a cut-off arguments. Roughly
speaking, an accurate uniformly boundedness estimates on
$(C)_c$ (Cerami) sequences of the associate energy functional $\Phi_\vr $
enables us to introduce a new functional $\widetilde{\Phi}_\vr $ by
virtue of the cut-off technique so that $\widetilde{\Phi}_\vr $ has
the same least energy solutions as ${\Phi}_\vr $ and can be dealt
with more easily under the assumption $\lm\in(0,\lm_0]$.

An outline of this paper is as follows: In section \ref{TVF} we treat
the linking argument which gives us a min-max scheme. In section
\ref{Preliminary results}, we study the limit equation and introduce the
cut-off arguments. Lastly, in section \ref{PMT},
the combination of the results in section \ref{TVF},
\ref{Preliminary results} proves the Theorem \ref{main theorem}.

%\textbf{Notation}: Throughout the paper we shall denote by $C,C_1,C_2$ etc. various
%positive constants which may vary from lines to lines and are not essential
%to the problem.

\section{The variational framework}\label{TVF}

\subsection{The functional setting and notations}

In the sequel, by $|\cdot|_q$ we denote the usual $L^q$-norm, and
$(\cdot,\cdot)_2$ the usual $L^2$-inner product. Let
$H_\om=i\al\cdot\nabla-a\bt+\om$ denote the self-adjoint operator on
$L^2\equiv L^2(\R^3,\C^4)$ with domain
$\cd(H_\om)=H^1\equiv H^1(\R^3,\C^4)$. It is
well know that $\sigma(H_\om)=\sigma_c(H_\om)=\R\setminus(-a+\om,a+\om)$
where $\sa(\cdot)$ and $\sa_c(\cdot)$ denote the spectrum and
the continuous spectrum. For $\om\in(-a,a)$, the space $L^2$ possesses the
orthogonal decomposition:
\begin{\equ}\label{l2dec}
L^2=L^+\oplus L^-,\ \ \ \ u=u^++u^-
\end{\equ}
so that $H_\om$ is positive definite (resp. negative definite) in
$L^+$ (resp. $L^-$). Let $E:=\mathcal{D}(\jdz{H_\om}^{1/2})=H^{1/2}$
be equipped with the inner product
$$
\inp{u}{v}=\Re(\jdz{H_\om}^{1/2}u,\jdz{H_\om}^{1/2}v)_2
$$
and the induced norm $\norm{u}=\inp{u}{u}^{1/2}$, where $\jdz{H_\om}$
and $\jdz{H_\om}^{1/2}$ denote respectively the absolute value of
$H_\om$ and the square root of $|H_\om|$. Since
$\sigma(H_\om)=\mathbb{R}\setminus(-a+\om,a+\om)$, one has
\begin{\equ}\label{l2ineq}
(a-|\om|)|u|_2^2\leq\norm{u}^2 \quad  \mathrm{for\ all\ }u\in E.
\end{\equ}
Note that this norm is equivalent to the usual $H^{1/2}$-norm, hence
$E$ embeds continuously into $L^q$ for all $q\in[2,3]$ and compactly
into $L_{loc}^q$ for all $q\in[1,3)$. It is clear that $E$ possesses
the following decomposition
\begin{\equ}\label{Edec}
E=E^+\oplus E^-\ \ \mathrm{with\ \ }E^{\pm}=E\cap L^{\pm},
\end{\equ}
orthogonal with respect to both $(\cdot,\cdot)_2$ and
$\inp{\cdot}{\cdot}$ inner products. This decomposition induces also
a natural decomposition of $L^p$, hence there is $d_p>0$ such that
\begin{\equ}\label{lpdec}
d_p\jdz{u^\pm}_p^p\leq\jdz{u}_p^p\ \ \mathrm{for\ all\ }u\in E.
\end{\equ}

Let $H^1(\R^3,\R)$ be equipped with the equivalent norm
$$
\norm{v}_{H^1}=\bigg (\int \jdz{\nabla v}^2+M v^2dx \bigg)^{1/2}
\quad \forall v\in H^1(\R^3,\R).
$$
Then (\ref{D3}) can be reduced to a single equation with a non-local
term. Actually, for any $v\in H^1$,
\begin{\equ}\label{R1}
\begin{aligned}
\jdz{4\pi \lm \int (\bt u)u \cdot v \, dx}&\leq
\bigg( 4\pi \lm \int |u|^2 |v| \, dx   \bigg)\\
 &\leq 4\pi \lm  |u|^2_{12/5} |v|_{6}\\
 &\leq 4\pi \lm S^{-1/2} |u|^2_{12/5}\|v\|_{H^1},
\end{aligned}
\end{\equ}
where $S$ is the Sobolev embedding constant: $S|v|^2_6\leq
\|v\|^2_{H^1}$ for all $v\in H^1$. Hence there
exists a unique $V_u\in H^1$
such that
\begin{\equ}\label{solution of poisson}
\int \nabla V_u \cdot \nabla z + M \cdot V_u z\, dx
=4\pi \lm \int (\bt u)u \cdot z\,dx
\end{\equ}
for all $z\in H^1$. It follows that $V_u$
satisfies the Schr\"{o}dinger type equation
\begin{\equ}\label{STE}
-\Delta V_u + M \cdot V_u=4\pi \lm (\bt u)u
\end{\equ}
and there holds
\begin{\equ}\label{juan ji}
V_u(x)=\lm \int_{\R^3}
\frac{[(\bt u)u](y)}{\jdz{x-y}}e^{-M|x-y|}\,dy.
\end{\equ}
Substituting $V_u$ in (\ref{D3}), we are led
to the equation
\begin{\equ}\label{D22}
H_\om u - \lm V_u \bt u  = W_\vr (x)g(|u|)u.
\end{\equ}

On $E$ we define the functional
$$
\Phi_\vr (u)=\frac{1}{2}\big( \|u^+\|^2-\|u^-\|^2 \big)
-\Gamma_\lm (u)-\Psi_\vr (u)
$$
for $u=u^++u^-$, where
$$
\Gamma_\lm (u)=\frac{\lm}{4}\int V_u\cdot (\bt u)u \, dx
=\frac{\lm^2}{4}\iint \frac{ [(\bt u)u](x) [(\bt u)u](y)}{|x-y|}
e^{-M|x-y|} \, dydx
$$
and
$$
\Psi_\vr (u)=\int W_\vr (x)G(\jdz{u})dx.
$$

\subsection{Technical results}

In this subsection, we shall introduce some lemmas related
to the functional $\Phi_\vr$.

\begin{Lem}\label{variation}
Under the hypotheses $(P_0)$, $(G_1)$-$(G_2)$,  one has
$\Phi_\vr \in C^2(E,\R)$ and any critical point
of $\Phi_\vr $ is a solution of \eqref{D3}.
\end{Lem}

\begin{proof}
Clearly, $\Psi_\vr\in C^2(E,\R)$. It remains to check that
$\Gamma_\lm\in C^2(E,\R)$. It suffices to show that, for any $u,v\in
E$,
\begin{\equ}\label{estimates of Gamma-eps0}
\jdz{\Gamma_\lm (u)}\leq C_1 \lm^2 \|u\|^4,
\end{\equ}
\begin{\equ}\label{estimates of Gamma-eps1}
\jdz{\Gamma_\lm '(u)v}\leq C_2 \lm^2 \|u\|^3 \|v\|,
\end{\equ}
\begin{\equ}\label{estimates of Gamma-eps2}
\jdz{\Gamma_\lm ''(u)[v,v]}\leq C_3 \lm^2 \|u\|^2 \|v\|^2.
\end{\equ}

Observe that one has, by using $V_u$ as a test function in
\eqref{STE},
\begin{\equ}\label{Ajsinequ}
|V_u|_6 \leq S^{-1/2}\|V_u\|_{H^1}
\leq C_1 \lm \|u\|^2.
\end{\equ}
This, together with the H\"older inequality (with $r=6, r'=6/5$),
implies \eqref{estimates of Gamma-eps0}. Note that $\Gamma'_\lm(u)v
= \frac{d}{dt}\Gamma_\lm(u+tv)\big|_{t=0}$, so
\begin{\equ}\label{R2}
\aligned \Gamma'_\lm(u)v =
&\frac{\lm^2}{2} \Re \iint\frac{e^{-M|x-y|}}{|x-y|}\Big(
[(\bt u)u](x) [(\bt u)v](y)   \\
& + [(\bt u)u](y) [(\bt u)v](x)
\Big)dydx\\
= &\, \lm \int V_u\cdot \Re (\bt u)v \,dx
\endaligned
\end{\equ}
which, together with the H\"older inequality and \eqref{Ajsinequ},
shows \eqref{estimates of Gamma-eps1}. Similarly,
\[\aligned
\Gamma^{''}_\lm(u)[v,v] = &\,2\,\lm^2
\iint\frac{e^{-M|x-y|}}{|x-y|}\Big(
\Re [(\bt u)v](x) \Re [(\bt u)v](y) \Big)dxdy \\
&\,\qquad
+ \lm \Re \int
V_u \cdot (\bt v)v  \, ,
\endaligned
\]
and one gets
\eqref{estimates of Gamma-eps2}.

Now it is a standard to verify that critical points of $\Phi_\vr$
are solutions of \eqref{D3}.
\end{proof}

\medskip

We show further the following:

\begin{Prop}\label{Gamma-eps nonnegative}
$\Gamma_\lm $ is non-negative and
weakly sequentially lower semi-continuous.
Moreover, $\Gamma_\lm$ vanishes only when
$(\bt u)u=0$ a.e. in $\R^3$.
\end{Prop}

\begin{proof}
Recall that for every $u\in E$, $V_u$ solves (in the weak sense)
$$
-\Delta V_u +M V_u = 4 \pi \lm (\bt u)u.
$$
Then a standard maximum principle argument shows that
\begin{\equ}\label{sere ineq}
\Big( V_u \cdot (\bt u)u
\Big) (x) \geq0, \quad {\rm a.e.\ on\ }\R^3.
\end{\equ}
Hence (see \eqref{juan ji})
$$
\Gamma_\lm (u)=\frac{\lm }{4}\int V_u\cdot
(\bt u)udx\geq0.
$$
Furthermore, suppose $u_n\rightharpoonup u$ in $E$,
then $u_n\rightarrow u$ a.e..
Therefore (\ref{sere ineq}) and Fatou's lemma yield
$$
\Ga_\lm (u)\leq\liminf_{n\to\infty}\Ga_\lm (u_n)
$$
as claimed.
\end{proof}

Set, for $r>0$, $B_r=\{u\in E:\norm{u}\leq r\}$, and for $e\in E^+$
$$E_e:=E^-\oplus\mathbb{R}^+e$$
with $\mathbb{R}^+=[0,+\infty)$. In virtue of the assumptions
$(G_1)$-$(G_2)$, for any $\de>0$, there exist $r_\de>0, c_\de>0$ and $c_\de'>0$ such that
\begin{\equ}\label{g estimates}
\left\{
\begin{aligned}
&g(s)<\de\ \ \mathrm{for\ all\ }0\leq s \leq r_\de;\\
&G(s)\geq c_\de\,s^\theta-\de\,s^2\ \ \mathrm{for\ all\ } s\geq 0;\\
&G(s)\leq \de\,s^2+c_\de'\,s^p\ \ \mathrm{for\ all\ } s\geq 0
\end{aligned}\right.
\end{\equ}
and
\begin{\equ}\label{R9}
\widehat{G}(s):=\frac12g(s)s^2-G(s)\geq\frac{\theta-2}{2\theta}g(s)s^2
\geq\frac{\theta-2}{2}G(s)\geq c_\theta s^\sigma
\end{\equ}
for all $s\geq 0$, where $c_\theta=c_0(\theta-2)/2$.

\begin{Lem}\label{max Phi-eps<C}
For all $\vr \in(0,1]$, $\Phi_\vr $ possess the linking structure:
\begin{itemize}
\item[$1)$] There are $r>0$ and $\tau>0$, both independent of $\vr $, such
that $\Phi_\vr |_{B_r^+}\geq0$ and $\Phi_\vr |_{S_r^+}\geq\tau$, where
\[
B_r^+=B_r\cap E^+=\{u\in E^+:\|u\|\leq r\},
\]
\[
S_r^+=\pa B_r^+=\{u\in E^+:\|u\|= r\}.
\]
\item[$2)$] For any $e\in E^+\setminus\{0\}$, there exist $R=R_e>0$ and
$C=C_e>0$, both independent of $\vr $, such that, for all $\vr >0$,
there hold $\Phi_\vr (u)<0$ for all $u\in E_e\setminus B_R$ and
$\max\Phi_\vr (E_e)\leq C$.
\end{itemize}
\end{Lem}

\begin{proof}
Recall that $\jdz{u}_p^p\leq C_p\norm{u}^p$ for all $u\in E$ by
Sobolev embedding theorem. 1) follows easily because, for $u\in E^+$
and $\de>0$ small enough
\[
\begin{aligned}
\Phi_\vr (u)&=\frac{1}{2}\norm{u}^2
 -\Gamma_\lm (u)-\Psi_\vr (u)\\
 &\geq\frac{1}{2}\norm{u}^2
 -C_1 \lm ^2 \|u\|^4
 -{|W|_\infty}\big(\de\jdz{u}_2^2+c_\de'\jdz{u}_p^p\big)\\
\end{aligned}
\]
with $C_1$, $C_p$ independent of $u$ and $p>2$ (see \eqref{estimates
of Gamma-eps0} and \eqref{g estimates}).

For checking 2), take $e\in E^+\setminus\{0\}$. In virtue of
(\ref{lpdec}) and \eqref{g estimates}, one gets, for $u=se+v\in
E_e$,
\begin{\equ}\label{linking ineq}
\begin{split}
\Phi_\vr (u)=&\,\frac{1}{2}\norm{se}^2-\frac{1}{2}\norm{v}^2
-\Gamma_\lm (u)-\Psi_\vr (u)\\
\leq&\, \frac{1}{2}s^2\norm{e}^2-\frac{1}{2}\norm{v}^2 - c_\de
d_\theta \inf W\cdot s^\theta\jdz{e}_\theta^\theta
\end{split}
\end{\equ}
proving the conclusion.
\end{proof}

Recall that a sequence $\{u_n\}\subset E$ is called to be a
$(PS)_c$-sequence for functional $\Phi\in C^1(E,\R)$
if $\Phi(u_n)\to c$ and $\Phi'(u_n)\to 0$, and is called to be
$(C)_c$-sequence for $\Phi$ if $\Phi(u_n)\to c$ and
$(1+\|u_n\|)\Phi'(u_n)\to 0$. It is clear
that if $\{u_n\}$ is a $(PS)_c$-sequence with $\{\|u_n\|\}$ bounded
then it is also a $(C)_c$-sequence. Below we are going to study
$(C)_c$-sequences for $\Phi_\vr$ but firstly we observe the
following

\begin{Lem}\label{A(u) bdd}
For all $u\in E$, we have
\[
\Bigg| \frac{V_u}{\lm \|u\|} \Bigg|_6 \leq
C |u|_\sa,
\]
where $\sa>0$ is the constant in $(G_2)$ and $C>0$ is
depending only on the embedding $H^1(\R^3,\R)\hookrightarrow L^6$
and $E\hookrightarrow L^q$ for
$\frac{1}{\sa}+\frac{1}{q}+\frac{1}{6}=1$.
\end{Lem}

\begin{proof}
Notice that
$V_u$ satisfies the equation
$$
-\Delta V_u + M V_u=4\pi \lm (\bt u)u,
$$
hence, using $V_u$ as a test function,
$$
\|V_u\|_{H^1}^2 \leq 4\pi \lm \int |V_u|\cdot |u|^2
$$
By H\"older's inequality
\[
\begin{aligned}
\|V_u\|_{H^1}^2 &\leq 4\pi \lm |V_u|_6 |u|_\sa |u|_q    \\
&\leq 4\pi \lm \tilde{C} \|V_u\|_{H^1} \cdot
\|u\|  \cdot |u|_\sa   \, .
\end{aligned}
\]
And then we infer
$$
\norm{\frac{V_u}{\lm\|u\|}}_{H^1}
\leq C |u|_\sa,
$$
which yields the conclusion.
\end{proof}

We now turn to an estimate on boundedness of $(C)_c$-sequences which
is the key ingredient in the sequel. Recall that, by $(G_1)$, there
exist $r_1>0$ and $a_1>0$ such that
\begin{\equ}\label{R8}
g(s)\leq \frac{a-|\omega|}{2\,|W|_\infty} \quad \text{for all $s\leq
r_1$},
\end{\equ}
and, for $s\geq r_1$,  $g(s)\leq a_1 s^{p-2}$, so
$g(s)^{\sigma_0-1}\leq a_2 s^2$ with
$$
\sigma_0:=\frac{p}{p-2}>3
$$
which, jointly with $(G_2)$, yields (see \eqref{R9})
\begin{\equ}\label{g-sigma0 estimate}
g(s)^{\sigma_0}\leq a_2 g(s)s^2\leq a_3 \widehat{G}(s)\quad
\text{for all $s\geq r_1$}.
\end{\equ}

\begin{Lem}\label{PScseq eps estimate}
Assume $(P_0)$, $(G_1)$-$(G_2)$ and $\lm>0$,
for every pair of constants $c_1,c_2>0$,
there exists a constant $\Lam>0$, depending only on
$c_1,c_2,\lm$, such that for any $u\in E$ with
\begin{\equ}\label{e1}
|\Phi_\vr(u)|\leq c_1
\quad {\rm and} \quad
\|u\|\cdot \|\Phi_\vr'(u)\| \leq c_2,
\end{\equ}
we have
\[
\|u\|\leq \Lam.
\]
Furthermore, $\Lam$ is a increasing function with respect to $\lm>0$.
\end{Lem}

Lemma \ref{PScseq eps estimate} has a immediate consequence
which implies the boundness of a $(C)_c$-sequence:

\begin{Cor}\label{PScseq2}
Consider $\vr\in(0,1]$, and $\{u_n^\vr\}$ is the corresponding
$(C)_{c_\vr}$-sequence for $\Phi_\vr$. If there exists $C>0$
such that $|c_\vr|\leq C$ for all $\vr$, then we have (up to a
subsequence if necessary)
\[
\|u_n^\vr\|\leq \Lam
\]
where $\Lam$ is found in Lemma \ref{PScseq eps estimate}
depends on $\lm$ and the pair $c_1=C$ and $c_2=1$.
\end{Cor}

\begin{proof}[Proof of Lemma \ref{PScseq eps estimate}]
Take $u\in E$ such that \eqref{e1} is satisfied. Without loss of
generality we may assume that $\|u\|\geq1$. The form of
$\Phi_\vr$ and the representation \eqref{R2} \,
($\Ga'_\lm(u)u=4\Ga_\lm(u)$) \, implies that
\begin{\equ}\label{R3}
\begin{aligned}
c_1+c_2\geq\Phi_\vr (u) -\frac{1}{2}\Phi_\vr '(u)u
=\Gamma_\vr (u)
 +\int W_\vr (x)\widehat{G}(|u|)\\
\end{aligned}
\end{\equ}
and
\begin{\equ}\label{R4}
\aligned
c_2\geq&\,\Phi'_\vr(u)(u^{+}-u^{-})\\
=&\,\|u\|^2-
\Gamma'_\lm(u)(u^{+}-u^{-})\\
&\, -\Re\int W_\vr(x)
 g(|u|)u \cdot (u^+-u^-).
\endaligned
\end{\equ}
By Lemma \ref{Gamma-eps nonnegative}, \eqref{R9} and \eqref{R3},
$|u|_\sa\leq C_1$, where $C_1$ depends only
 on $c_1,c_2$. It follows from \eqref{R4} that
\[
\|u\|^2 \leq  c_2+
\Gamma'_\lm(u )(u^+
 -u^-)
 +\Re\int W_\vr (x)
 g(|u|)u \cdot (u^+-u^-).
\]
This, together with \eqref{R8} and \eqref{l2ineq}, shows
\begin{\equ}\label{R5}
\frac{1}{2}\|u\|^2  \leq c_2+
\Gamma'_\lm(u)(u^+ - u^-)
 +\Re\int_{|u|\geq r_1} W_\vr (x)
 g(|u|)u \cdot (u^+ -u^-).
\end{\equ}

Recall that $(G_1)$ and $(G_2)$ imply $2<\sa\leq p$. Setting
$t=\frac{p\sigma}{2\sigma-p}$, one sees
\[
2 < t < p, \quad \frac{1}{\sigma_0}+\frac{1}{\sigma}+\frac{1}{t}=1.
\]
By H\"older inequality, the fact $\Gamma_\lm(u)\geq 0$,
\eqref{g-sigma0 estimate}, \eqref{R3} and the embedding of $E$
to $L^t$, we have
\begin{\equ}\label{R6}
\aligned
&\, \int_{|u|\geq r_1} W_\vr(x) \, g(\jdz{u
})|u|\cdot |u^+-u^-|     \\
\leq &\, |W|_\infty\Big(\int_{|u|\geq r_1}
g(|u|)^{\sa_0}\Big)^{1/\sa_0} \Big(\int
|u|^\sa\Big)^{1/\sa} \Big(|u^+-
u^-|^t\Big)^{1/t}\\
\leq &\, C_2\|u\|
\endaligned
\end{\equ}
with $C_2>0$ depends only on $c_1,c_2$.

Let $q=\frac{6\sigma}{5\sigma-6}$. Then $2<q<3$ and
$\frac{1}{\sigma}+\frac{1}{q}+\frac{1}{6}=1$. Set
\[
\zeta=\left\{\aligned & 0 \, & \text{if $q=\sigma$};\\
& \frac{2(\sigma-q)}{q(\sigma-2)}  \, & \text{if $q<\sigma$};\\
& \frac{3(q-\sigma)}{q(3-\sigma)}  \, & \text{if $q>\sigma$};
\endaligned\right.
\]
we deduce that $\zeta<1$ and
$$
|u|_q \leq \left\{\aligned & \jdz{u}_2^{\zeta}
\cdot\jdz{u}_\sigma^{1-\zeta} \, & \mathrm{if}\ 2<q\leq \sigma\\
& \jdz{u}_3^{\zeta}\cdot\jdz{u}_\sigma^{1-\zeta} \, & \mathrm{if}\
\sigma< q<3.
\endaligned\right.
$$
By virtue of the H\"older inequality, Lemma \ref{Gamma-eps
nonnegative} and the
embedding of $E$ to $L^2$ and $L^3$, we obtain
\[
\begin{aligned}
& \Bigg| \lm \Re\int V_u\cdot (\bt u) (u^+-u^-) \Bigg|\\
\leq &\lm \|u\|  \Bigg|  \Re\int \frac{V_u}{\|u\|}
(\bt u) \cdot (u^+-u^-)  \Bigg|  \\
\leq&\lm^2 \|u\| \Bigg| \frac{V_u}{\lm \|u\|} \Bigg|_6
|u|_\sa \cdot |u^+-u^-|_q\\
\leq &\, \lm^2 C_3\|u\| \cdot |u|_q\,\leq \,\lm^2
C_4\|u\|^{1+\zeta}
\end{aligned}
\]
with $C_4>0$ depends only on the embedding $E\hookrightarrow L^q$.
This, together with the representation of \eqref{R2}, implies that
\begin{\equ}\label{R7}
|\Ga'_\lm(u)(u^+-u^-)|\leq
\lm^2  C_4\|u\|^{1+\zeta}  \, .
\end{\equ}

Now the combination of \eqref{R5}, \eqref{R6} and \eqref{R7} shows
that
%\begin{\equ}\label{x1}
%\frac{M_1}{\norm{u_n^\vr }^{1-\zeta}}+\frac{ M_2}{\norm{u_n^\vr
%}}\geq \, 1
%\end{\equ}
\begin{\equ}\label{x1}
\|u\|^2\leq M_1 \|u\| + \lm^2 M_2 \|u\|^{1+\zeta}
\end{\equ}
with $M_1$ and $M_2$ dependent only on the constants $c_1, c_2$.
Therefore,
either $\|u\|\leq 1$ or there is $\Lam\geq 1$ dependents
only on $c_1, c_2, \lm$ such that
\[
\|u\|\leq \Lam
\]
as desired. Moreover, \eqref{x1} implies $\Lam$ is increasing in $\lm$.
\end{proof}

Finally, for later aims we define the operator $\mathcal{V}
:E\to H^1(\R^3,\R)$ by $\mathcal{V}(u)=V_u$. We have
\begin{Lem}\label{lemma of Aj}
\begin{itemize}
\item[$(1)$] $\mathcal{V}$ maps bounded sets into bounded sets;
\item[$(2)$] $\mathcal{V}$ is continuous;
\end{itemize}
\end{Lem}

\begin{proof}
Clearly, (1) is a straight consequence of (\ref{Ajsinequ}). (2)
follows easily because, for $u,v\in E$, one sees that $V_u-V_v$ satisfies
$$
-\Delta(V_u-V_v)+M(V_u-V_v)=
4\pi \lm [(\bt u)u-(\bt v)v].
$$
Hence
\[
\begin{aligned}
\|V_u-V_v\|_{H^1}&\leq \lm C
\big| (\bt u)u - (\bt v)v \big|_{6/5}\\
 &\leq \lm C \Big(\jdz{u-v}_{12/5}\jdz{u}_{12/5}
 +\jdz{u-v}_{12/5}\jdz{v}_{12/5}\Big)\\
 &\leq \lm \tilde{C} \big(
 \|u-v\|\cdot\|u\|+\|u-v\|\cdot\|v\|
 \big).
\end{aligned}
\]
and this implies the desired conclusion.
\end{proof}

\section{Preliminary results}\label{Preliminary results}

We are interested in describing the concentration phenomena
of the least energy solutions to the semi-classical model
\eqref{D3}. Throughout this section we will collect properties
of the energy functionals of the Dirac-Klein-Gordon systems
(including the estimates of the least energy). Instead of
dealing directly with the nonlocal term $\Ga_\lm$, it seems
simpler to consider a modified problem
(see subsection \ref{modified problem}). For reasons that
will be apparent later, we treat our model in the case $\lm$
is not chosen large, that is $\lm\in(0,\lm_0]$ for some
$\lm_0>0$ will be chosen later on.

\subsection{The limit equation}\label{subsec1}

In order to prove our main result, we
will make use of the limit equation. For any $\mu>0$, consider the
equation
\begin{\equ}\label{limt equ mu}
\left\{
\aligned
&i\alpha\cdot\nabla u - a\bt u + \om u - \lm V \bt u = \mu g(\jdz{u})u,\\
&-\De V +M \cdot V=4\pi \lm (\bt u)u.
\endaligned \right.
\end{\equ}
Its solutions are critical points of the functional
\[
\begin{aligned}
\mst_\mu(u)&:=\frac{1}{2}\big(
\|u^+\|^2-\|u^-\|^2  \big)  - \Ga_\lm(u)
 -{\mu}\int G(\jdz{u})\\
 &\, =\frac{1}{2}\big(
 \|u^+\|^2-\|u^-\|^2  \big) - \Ga_\lm(u)
 -\msg_\mu(u)
\end{aligned}
\]
defined for $u=u^++u^-\in E=E^+\oplus E^-$. Denote the critical
set and the least energy  of
$\mst_\mu$ as follows
\[
\aligned
&\msk_\mu:=\{u\in E:\ \mst_\mu'(u)=0\},\\
&\ga_\mu:=\inf\{\mst_\mu(u):\
 u\in\msk_\mu\setminus\{0\}\}.
\endaligned
\]

In order to find critical points of $\mst_\mu$, we will
use the following abstract theorem
which is taken from \cite{Ding1, Ding2}.

Let $E$ be a Banach space with direct sum decomposition
$E = X \op Y$, $u = x + y$
and corresponding projections $P_X , P_Y$
onto $X, Y$, respectively. For a functional
$\Phi\in C^1(E,\R)$ we write $\Phi_a = \{
u \in E : \, \Phi(u)\geq a \}$.

Now we assume that $X$ is separable and reflexive,
and we fix a countable dense subset $\mss \subset X^*$.
For each $s \in \mss$ there is a semi-norm on E defined by
\[
p_s: E \to \R , \quad
p_s(u)= |s (x)| + \| y \|  \quad
{\rm for}\  u = x + y \in X \op Y.
\]
We denote by $\ct_\mss$
the induced topology. Let $w^*$
denote the weak*-topology on $E$.
Suppose:
\begin{itemize}
\item[$(\Phi_0)$] There exists $\xi > 0$ such that
$\|u\| <\xi \|P_Y u\|$ for all $u\in\Phi_0$.

\item[$(\Phi_1)$] For any $c\in\R$, $\Phi_c$ is $\ct_\mss$-closed, and
$\Phi':(\Phi_c,\ct_\mss)\to(E^*,w^*)$ is continuous.

\item[$(\Phi_2)$] There exists $\rho>0$ with
$\kappa:= \inf \Phi( S_\rho Y ) > 0$ where $S_\rho Y:=\{ u\in Y : \, \|u\| = \rho \}$.
\end{itemize}
The following theorem is a special case of \cite{Ding1} Theorem 3.4
(see also \cite{Ding2} Theorem 4.3).

\begin{Thm}\label{CPT}
Let $(\Phi_0)-(\Phi_2)$ be satisfied and suppose there are
$R>\rho>0$ and $e\in Y$ with $\|e\| = 1$ such that $\sup\Phi(\pa Q)\leq \ka$
where $Q = \{ u = x + te :\, x \in X, t \geq 0 , \|u\|<R \}$ .Then $\Phi$ has a
$(C)_c$-sequence with $\ka\leq c\leq\sup\Phi( Q )$.
\end{Thm}

The following lemma is useful to verify $(\Phi_1)$
(see \cite{Ding1} or \cite{Ding2}).

\begin{Lem}\label{Phi 1}
Suppose $\Phi\in C^1(E,\R)$ is of the form
\[
\Phi(u)=\frac{1}{2} \big(
\|y\|^2-\|x\|^2  \big)-\Psi(u) \quad
{\rm for}\  u= x + y \in E = X \op Y
\]
such that
\begin{itemize}
\item[(i)] $\Psi\in C^1(E,\R)$ is bounded from below;

\item[(ii)] $\Psi:(E,\ct_w)\to\R$ is sequentially lower semi-continuous, that is,
$u_n\rightharpoonup u$ in $E$ implies $\Psi(u) \leq \liminf \Psi(u_n)$;

\item[(iii)] $\Psi':(E,\ct_w)\to(E^*,w^*)$ is sequentially continuous.

\item[(iv)] $\nu: E\to\R$, $\nu(u)=\|u\|^2$, is $C^1$
and $\nu\,':(E,\ct_w)\to(E^*,w^*)$  is sequentially
continuous.
\end{itemize}
Then $\Phi$ satisfies $(\Phi_1)$.
\end{Lem}

Next, we present the existence result for the limit equation
\eqref{limt equ mu}.

\begin{Lem}\label{erfle}
Let $\lm$ be a positive constant, for each $\mu>0$, we have
\begin{enumerate}
\item $\msk_\mu\not=\emptyset$ and $\ga_\mu>0$,

\item $\ga_\mu$ is attained.
\end{enumerate}
\end{Lem}

\begin{proof}
Invoking Lemma \ref{Gamma-eps nonnegative}, we see
that $(\Phi_0)$ is satisfied.
With $X=E^-$ and $Y=E^+$ the condition $(\Phi_0)$ holds
by Lemma \ref{Gamma-eps nonnegative} and Lemma \ref{Phi 1}.
Together with the linking structure (see Lemma \ref{max Phi-eps<C})
we have all the assumptions of Theorem \ref{CPT} verified. Therefore,
there exists a sequence $\{u_m\}$ satisfying $\mst_\mu(u_m)\to c > 0$
and $(1+\|u_m\|)\mst_\mu'(u_m)\to0$ as $m\to\infty$. Using the same
arguments in proving Lemma \ref{PScseq eps estimate}, we get $\{u_m\}$
is bounded. Now by the classical concentration compactness principle
(cf. \cite{Lions}) and the translation-invariance of $\mst_\mu$,
we infer there is $u\neq0$ such that $\mst_\mu'(u)=0$.

If $u\in\msk_\mu$, one has
\begin{\equ}\label{XX1}
\mst_\mu(u)=\mst_\mu(u)-\frac{1}{2}\mst_\mu'(u)u
=\Ga_\lm(u)+\mu\int \widehat{G}(u)\geq0.
\end{\equ}
For proving $\ga_\mu>0$, assume by contradiction that $\ga_\mu=0$.
Let $u_j\in\msk_\mu\setminus\{0\}$ such that $\mst_\mu(u_j)\to0$.
It is obvious that $\{u_j\}$ is bounded. Furthermore, by \eqref{R9}
and \eqref{XX1}, we deduce $u_j\to0$ in $L^\sa$ as $j\to \infty$.
On the other hand, by noting that
$0=\mst_\mu'(u_j)(u_j^+-u_j^-)$, \eqref{lpdec} and Lemma \ref{A(u) bdd}
imply
\[
\aligned
\|u_j\|^2&=\Ga_\lm'(u_j)(u_j^+-u_j^-)+\mu\int g(u_j)(u_j^+-u_j^-)\\
 & \leq \lm^2 C_1 \|u_j\|^3 \cdot |u_j|_\sa +\mu\int g(u_j)(u_j^+-u_j^-) \, .
\endaligned
\]
By \eqref{g-sigma0 estimate} and H\"older's inequality, one sees
\[
\aligned
\frac{1}{2}\|u_j\|^2 &\leq \lm^2 C_1 \|u_j\|^3 \cdot |u_j|_\sa
+ C_2 \mu \bigg( \int g(|u_j|)^{\sa_0} \bigg)^{1/\sa_0} |u_j|_p^2 \\
&\leq \lm^2 C_1 \|u_j\|^3 \cdot |u_j|_\sa
+ C_3 \mu \big( \mst_\mu(u_j) \big)^{1/\sa_0} \|u_j\|^2 \,  .
\endaligned
\]
Hence $\frac{1}{2}\leq o(1)+o(1)$, a contradiction.

Lastly, again, by using the concentration compactness principle,
we check easily that $\ga_\mu$ is attained, ending the proof.
\end{proof}

\subsection{A modification for the nonlocal term}\label{modified problem}

We find our current research is more delicate, since
the solutions we look for are at the least energy level and
$\Ga_\lm$ is not convex on $E$ (even for $u$ with $\|u\|$ large).
By cutting off  the nonlocal terms, we are able to find a
critical point via an appropriate min-max scheme.
The critical point will eventually
be shown to be a least energy solution to our model.

Next we introduce the modified problem by
choosing a cut-off function $\eta:\R\to\R$ such that
$\msf_\lm(u):=\eta(\|u\|^2)\Ga_\lm(u)$ vanishes for $\|u\|$
large.

By virtue of $(P_0)$, set $b=\inf W(x)>0$, let us first consider
the autonomous systems for $\mu\geq b$
\[
\left\{
\aligned
&i\alpha\cdot\nabla u - a\bt u + \om u - \lm V \bt u = \mu g(\jdz{u})u,\\
&-\De V +M \cdot V=4\pi \lm (\bt u)u.
\endaligned \right.
\]
Following Lemma \ref{erfle}, $\ga_\mu>0$ (the least energy)
is attained. Now fix $\Lam>0$ to be the constant
(independent of $\vr>0$) found in Lemma \ref{PScseq eps estimate}
associated to $\lm>0$ and the pair of the constant $c_1=C_{e_0}$ and
$c_2=1$, where $C_{e_0}$ (independent of $\lm$ and $\mu$)
is the constant in Lemma \ref{max Phi-eps<C} with ${e_0}\in E^+\setminus\{0\}$
being fixed.

It is obvious that $\ga_\mu\leq C_{e_0}$.
%Now we can see $\Lam$ as
%a monotonous function with respect to $\lm>0$.
Denote $T=(\Lam+1)^2$
and choose $\eta:[0,+\infty)\to[0,1]$ be a smooth function with
$\eta(t)=1$ if $0\leq t \leq T$, $\eta(t)=0$ if $t\geq T+1$,
$\max |\eta'(t)|\leq 2$ and $\max|\eta''(t)|\leq 2$. Define
$\msf_\lm:E\to \R$ as $\msf_\lm(u)=\eta(\|u\|^2)\Ga_\lm(u)$.
Then we have $\msf_\lm\in C^2(E,\R)$ and $\msf_\lm$ vanishes for
all $u$ with $\|u\|\geq \sqrt{T+1}$.

Consider the modified functionals
\[
\widetilde{\mst}_\mu(u)=\frac{1}{2}\big(
\|u^+\|^2-\|u^-\|^2 \big) - \msf_\lm(u) - \msg_\mu(u),
\]
and
\[
\widetilde{\Phi}_\vr(u)=\frac{1}{2}\big(
\|u^+\|^2-\|u^-\|^2 \big) - \msf_\lm(u) - \Psi_\vr(u).
\]
By definition, $\widetilde{\mst}_\mu\big|_{B_T}=\mst_\mu$
and $\widetilde{\Phi}_\vr\big|_{B_T}=\Phi_\vr$ where
$B_T:=\{u\in E: \|u\|\leq \sqrt T\}$. And it's easy to see that
$0\leq \msf_\lm(u)\leq \Ga_\lm(u)$ and
\[
\big|\msf_\lm '(u)v \big| \leq \big| 2\eta'(\norm{u}^2)\Ga_\lm
(u)\inp{u}{v} \big| + \big| \Ga_\lm '(u)v \big|
\]
for $u,v\in E$.

Similarly to Lemma \ref{PScseq eps estimate}, we have
the following boundedness lemma (with $\Lam$ being
taken as above):

\begin{Lem}\label{bdd mod}
Assume $(G_1)$-$(G_2)$ and $(P_0)$. There exists $\lm_1>0$ such that,
for each $\lm\in(0,\lm_1]$, if $u\in E$ satisfies
\begin{\equ}\label{XX2}
0\leq \widetilde{\Phi}_\vr(u) \leq C_{e_0} \quad {\rm and} \quad
\|u\|\cdot \big\|\widetilde{\Phi}_\vr' \big\| \leq 1,
\end{\equ}
then we have $\|u\|\leq \Lam+1$, and consequently
$\widetilde{\Phi}_\vr(u)=\Phi_\vr(u)$.

In particular, replace $\widetilde{\Phi}_\vr$ with $\widetilde{\mst}_\mu$,
we have $\widetilde{\mst}_\mu$ shares the same ground state solution with $\mst_\mu$.
\end{Lem}

\begin{proof}
We repeat the arguments of Lemma \ref{PScseq
eps estimate}. Let $u$ satisfy \eqref{XX2}. If $\|u\|^2\geq T+1$
then $\msf_\lm(u)=0$ so, as proved in Lemma
\ref{PScseq eps estimate}, one changes \eqref{x1} by
$\|u\|^2\leq M_1 \|u\|$ and gets $\|u\|\leq \Lam$,
a contradiction.
Thus we assume that $\|u\|^2\leq T+1$. Then,
using \eqref{estimates of Gamma-eps0},  $ |\eta'(\|u\|^2)
\|u\|^2 \Ga_\lm(u)| \leq \lm^2 d_\lm^{(1)} $ (here and in the
following, by $d_\lm^{(j)}$ we denote positive constants depending only on
$\lm$ and $d_\lm^{(j)}$ is increasing with respect to $\lm$).
Similar to \eqref{R3},
\[
C_{e_0}+1\geq
\big(\eta(\|u\|^2)+2\eta'(\|u\|^2)\|u\|^2\big)\Ga_\lm(u)+
\int W_\vr(x)\widehat{G}(|u|)
\]
which yields
%\begin{\equ}\label{R9}
\[
C_{e_0}+1+\lm^2d_\lm^{(1)}> \eta(\|u\|^2)\Ga_\lm(u)+ \int
W_\vr(x)\widehat{G}(|u|),
\]
consequently $|u|_\sa\leq d_\lm^{(2)}$. Similarly to \eqref{R5} we
get that
\[
\aligned
%&\, o(1)+
\frac12\|u\|^2 \leq &\,\lm^2d_\lm^{(3)} +
 \eta(\|u\|^2)\Ga'_\lm(u)(u^+-u^-)\\
 &\,+\Re\int_{|u|\geq r_1} W_\vr (x)
 g(|u|)u \cdot\overline{u^+-u^-}
\endaligned
\]
which, together with \eqref{R6} and \eqref{R7}, implies either
$\|u\|\leq 1$ or as \eqref{x1}
\[
\|u\|^2\leq \lm^2 d_\lm^{(4)}+M_1\|u\|+M_2\|u\|^{1+\zeta},
\]
thus
\[
\|u\|\leq \lm^2d_\lm^{(5)}+\Lam.
\]
By monotonicity of $d_\lm^{(j)}$, we see that, for $\lm_1>0$ being suitably chosen,
let $\lm\in(0,\lm_1]$ then $\|u\|\leq \Lam+1$.
The proof is complete.
\end{proof}

\subsection{Estimates on the least energy}

Under Lemma \ref{bdd mod}, instead of study directly on $\Phi_\vr$ and
$\mst_\mu$, we turn to investigate the modified functionals, that is,
$\widetilde{\Phi}_\vr$ and $\widetilde{\mst}_\mu$ respectively. This
will give more information on the least energy level and more descriptions
on the min-max scheme.

Firstly, following the definitions of the modified functionals, an easy observation shows:

\begin{Prop}\label{PP1}
$\widetilde{\Phi}_\vr$ and $\widetilde{\mst}_\mu$
possess the linking structure proved in Lemma
\ref{max Phi-eps<C}, and the constants found in Lemma \ref{max Phi-eps<C}
are independent of the choice of $\widetilde{\Phi}_\vr$, $\Phi_\vr$,
$\widetilde{\mst}_\mu$ or
$\mst_\mu$, where $\mu\geq b$.
\end{Prop}

Now let us define (see \cite{Ding1,Szulkin})
\begin{\equ}\label{minmax level}
c_\vr :=\inf_{z\in E^+\setminus\{0\}}
\max_{u\in E_z}\widetilde{\Phi}_\vr (u)  \quad {\rm and} \quad
\tilde{\ga}_\mu:=\inf_{z\in E^+\setminus\{0\}}
\max_{u\in E_z}\widetilde{\mst}_\mu (u)
\end{\equ}
As a consequence of Proposition \ref{PP1} and Lemma \ref{bdd mod}
we have

\begin{Lem}\label{C independent of eps}
$c_\vr, \tilde{\ga}_\mu\in[\tau, C_{e_0}]$. Moreover, consider $\mu\geq b$,
if $c_\vr$ and $\tilde{\ga}_\mu$ are critical values for $\widetilde{\Phi}_\vr$
and $\widetilde{\mst}_\mu$, then they are also critical values for
$\Phi_\vr$ and $\mst_\mu$ respectively.
\end{Lem}

For a specific description, let us introduce the following notations:
Consider $\mu\geq b$, define
\[
\ci=\left\{
\aligned
&\widetilde{\Phi}_\vr \quad {\rm for\ the\ nonautonomous\ system},\\
&\widetilde{\mst}_\vr \quad {\rm for\ the\ autonomous\ system}.
\endaligned \right.
\]
Following Ackermann \cite{Ackermann} (also see
\cite{Ding2010,Ding2012,Ding2008}), for any fixed $u\in E^+$,
let $\varphi_u:E^-\to\mathbb{R}$
defined by $\varphi_u(v)=\ci(u+v)$. We have, for any
$v,w\in E^-$,
\[
\aligned
\varphi_u''(v)[w,w]
%&= -\norm{w}^2-\msf_\lm''(u+v)[w,w]
% -\msg_\mu''(u+v)[w,w]   \\
\leq  -\norm{w}^2-\msf_\lm''(u+v)[w,w]  \,  .
\endaligned
\]
At this point, a direct computation shows
\[
\aligned
& \msf_\lm (u+v)''[w,w]\\
=&\bkt{4\eta''(\norm{u+v}^2)\jdz{\inp{u+v}{w}}^2
+2\eta'(\norm{u+v}^2)\norm{w}^2}\Ga_\lm (u+v)\\
&+4\eta'(\norm{u+v}^2)\inp{u+v}{w}\Ga_\lm '(u+v)w
+\eta(\norm{u+v}^2)\Ga_\lm ''(u+v)[w,w]  \,  .
\endaligned
\]
Combining (\ref{estimates of Gamma-eps0})-(\ref{estimates of
Gamma-eps2}) yields
\[
\big|
\msf_\lm ''(u+v)[w,w] \big|
\leq \lm^2 d_\lm \|w\|^2 \leq \frac12 \|w\|^2
\]
for $\lm\leq\lm_2$, where $\lm_2$ is suitably chosen (here $d_\lm$ is
a positive constant depending monotonically only on $\lm$).
Hence, by setting $\lm_0=\min\{\lm_1,\lm_2\}$, for each $\lm\in(0,\lm_0]$
we deduce
\[
\varphi_u''(v)[w,w]
\leq  -\frac12 \|w\|  \,  .
\]
Additionally, we find
$$
\varphi_u(v)\leq\frac12 \big( \|u\|^2
 -\|v\|^2 \big).
$$
Therefore, there exists a unique $\xi:E^+\to E^-$ such that
$$
\ci(u+\xi(u))=\max_{v\in E^-}\ci(u+v).
$$
Here we used the expressions
\[
\xi(u)=\left\{
\aligned
&h_\vr(u) \quad {\rm defined\ for\ the\ nonautonomous\ system},\\
&\msj_\mu(u) \quad {\rm defined\ for\ the\ autonomous\ system}.
\endaligned \right.
\]
In the sequel, we fix $\lm$ in the interval $(0,\lm_0]$.
Next, setting $I_\vr, J_\mu : E^+\to\R$ by
\[
\aligned
&I_\vr(u)=\widetilde{\Phi}_\vr(u+h_\vr(u)),\\
&J_\mu(u)=\widetilde{\mst}_\mu(u+\msj_\mu(u)),
\endaligned
\]
and
\[
\aligned
&\msn_\vr=\{u\in E^+\setminus\{0\}:\
 I_\vr'(u)u=0\}, \\
&\msm_\mu=\{u\in E^+\setminus\{0\}:\
 J_\mu'(u)u=0\} \, .
\endaligned
\]
Denote by
\[
\cj(u)=\left\{
\aligned
&I_\vr(u) \quad {\rm for\ the\ nonautonomous\ system},\\
&J_\mu(u) \quad {\rm for\ the\ autonomous\ system},
\endaligned \right.
\]
and
\[
\cm=\left\{
\aligned
&\msn_\vr \quad {\rm for\ the\ nonautonomous\ system}\, , \\
&\msm_\mu \quad {\rm for\ the\ autonomous\ system}\, .
\endaligned \right.
\]
Plainly, critical points of $\cj$ and $\ci$ are in one-to-one
correspondence via the injective map $u\mapsto u+\xi(u)$ from
$E^+$ into $E$.

\begin{Lem}\label{unique t}
For any $u\in E^+\setminus\{0\}$, there is a unique $t=t(u)>0$ such
that $tu\in\cm$.
\end{Lem}

\begin{proof}
See \cite{Ackermann,Ding2008}.
\end{proof}

To give more information on the min-max levels defined in
\eqref{minmax level}, we set
\[
d=\left\{
\aligned
& c_\vr \quad {\rm for\ the\ nonautonomous\ system}\, ,\\
& \tilde{\ga}_\mu \quad {\rm for\ the\ autonomous\ system}\, .
\endaligned \right.
\]

\begin{Prop}\label{d-eps=c-eps}
There holds:
\begin{enumerate}
\item $d =\inf_{u\in\cm }\cj (u)$.

\item For $\mu\geq b$,
$\tilde{\ga}_\mu$ is the least energy for $\widetilde{\mst}_\mu$
and, by invoking Lemma \ref{bdd mod}, $\tilde{\ga}_\mu=\ga_\mu$.

\item Let $u\in\msm_\mu$ be such that $J_\mu(u)=\tilde{\ga}_\mu$
and set $E_u=E^-\op \R^+u$. Then
\[
\max_{w\in E_u}\widetilde{\mst}_\mu(w)=J_\mu(u).
\]

\item If $\mu_2>\mu_1\geq b$, then $\tilde{\ga}_{\mu_1}>\tilde{\ga}_{\mu_2}$.
\end{enumerate}
\end{Prop}

\begin{proof}
Denoting $\hat d =\inf_{u\in\cm }\cj (u)$,
given $e\in E^+$, if $u=v+se\in E_e$ with $\cj(u)=\max_{z\in E_e}\ci (z)$
then the restriction
$\ci |_{E_e}$ of $\ci $ on $E_e$
satisfies $(\ci |_{E_e})'(u)=0$ which implies
$v=\xi(se)$ and $\ci'(se)(se)=0$, i.e. $se\in\cm$.
Thus $\hat d \leq d $. While, on the other hand, if
$w\in\cm $ then $(\ci |_{E_w})'(w+\xi
(w))=0$, hence, $d \leq\max_{u\in E_w}\ci
(u)=\cj (w)$. Thus $\hat d \geq d $. It follow that $d=\hat d$.
Since it is standard to see that, for the autonomous system,
$\inf_{u\in\msm_\mu }J_\mu (u)$ characterize the least energy, we
infer that $\ga_\mu=\tilde{\ga}_\mu$. To prove 3, we note that
$u+\msj_\mu(u)\in E_u$ and
\[
J_\mu(u)=\widetilde{\mst}_\mu(u+\msj_\mu(u))\leq
\max_{w\in E_u} \widetilde{\mst}_\mu(w),
\]
moreover, since $u\in\msm_\mu$,
\[
\max_{w\in E_u}\widetilde{\mst}_\mu(w) \leq
\max_{s\geq0}\widetilde{\mst}_\mu(su+\msj_\mu(su))\leq
\max_{s\geq0}J_\mu(su)=J_\mu(u).
\]
Therefore, $\max_{w\in E_u}\widetilde{\mst}_\mu(w)=J_\mu(u)$.
Lastly to get 4, let $u_1$ be the ground state solution
for $\widetilde{\mst}_{\mu_1}$
and set $e=u_1^+$. Then
\[
\tilde{\ga}_{\mu_1}=\widetilde{\mst}_{\mu_1}(u_1)
=\max_{w\in E_e}\widetilde{\mst}_{\mu_1}(w).
\]
Suppose $u_2\in E_e$ be such that
$\widetilde{\mst}_{\mu_2}(u_2)=\max_{w\in E_e}\widetilde{\mst}_{\mu_2}(w)$.
We deduce that
\[
\aligned
\tilde{\ga}_{\mu_1}
&=\widetilde{\mst}_{\mu_1}(u_1)\geq
 \widetilde{\mst}_{\mu_1}(u_2)=\widetilde{\mst}_{\mu_2}(u_2)
 +(\mu_2-\mu_1)\int G(|u_2|)\\
&\geq \tilde{\ga}_{\mu_2}+(\mu_2-\mu_1)\int G(|u_2|).
\endaligned
\]
This ends the proof.
\end{proof}

\begin{Lem}\label{For any e  Te}
For any $e\in E^+\setminus\{0\}$,
there is $T_e>0$ independent the choice
of $\widetilde{\Phi}_\vr$ or
$\widetilde{\mst}_\mu$ such that
$t_e \leq T_e$ for $t_e >0$  satisfying
$t_e  e\in\cm $.
\end{Lem}

\begin{proof}
Since $\cj '(t_e  e)(t_e  e)=0$, one
get
$$
\ci (t_e  e+\xi (t_e  e))
=\max_{w\in E_e}\ci (w)\geq\tau.
$$
This,
together with Proposition \ref{PP1} (the linking structure), shows the
assertion.
\end{proof}

\subsection{Some auxiliary results}

Now using the notations introduced above, we are going to show some
auxiliary results that will make our arguments more transparent. First
of all, to describe the nonlinearities, we set
\[
\cn(u)=\left\{
\aligned
&\Psi_\vr(u) \quad {\rm for\ the\ nonautonomous\ system}\, ,\\
& \msg_\mu(u) \quad {\rm for\ the\ autonomous\ system}\, .
\endaligned \right.
\]

For any $u\in E^+$ and $v\in E^-$, setting $z=v-\xi(u)$ and
$l(t)=\ci(u+\xi(u)+tz)$, one has
$l(1)=\ci(u+v)$, $l(0)=\ci(u+\xi(u))$ and
$l'(0)=0$. Thus $l(1)-l(0)=\int_0^1(1-t)l''(t)dt$. This implies that
\[
\begin{aligned}
&\ci(u+v)-\ci(u+\xi(u))\\
=&\int_0^1(1-t)\ci''\bkt{u+\xi(u)+tz}[z,z]dt\\
=&-\int_0^1(1-t)\|z\|^2 dt
-\int_0^1(1-t)\Big[
\msf_\lm''(u+\xi(u)+tz)[z,z]\\
&+\cn''(u+\xi(u)+tz)[z,z] \Big] dt,
\end{aligned}
\]
and hence
\begin{\equ}\label{the equa 1}
\begin{split}
 &\int_0^1(1-t)\Big[
\msf_\lm''(u+\xi(u)+tz)[z,z]
+\cn''(u+\xi(u)+tz)[z,z] \Big] dt\\
 &+\frac{1}{2}\|z\|^2
 =\ci(u+\xi(u))-\ci(u+v).
\end{split}
\end{\equ}

\begin{Rem}\label{auto inequ X}
Recall that, for $\lm\in(0,\lm_0]$ being a positive constant,
there holds
\[
\big| \msf_\lm''(u+\xi(u)+tz)[z,z] \big| \leq \frac12\|z\|^2.
\]
From \eqref{the equa 1}, we deduce that, for the autonomous system,
\begin{\equ}\label{auto inequ}
\aligned
&\widetilde{\mst}_\mu(u+\msj_\mu(u))-\widetilde{\mst}_\mu(u+v) \\
\geq &\, \frac14 \|z\|^2
 + \int_0^1(1-t)\msg_\mu''(u+\xi(u)+tz)[z,z] dt \, .
\endaligned
\end{\equ}
\end{Rem}

Next we estimate the regularity of the critical points of
$\widetilde{\Phi}_\vr$.  Let $\msk_\vr :=\{u\in
E:\ \widetilde{\Phi}_\vr '(u)=0\}$ be the critical set of
$\widetilde{\Phi}_\vr $. It is easy to see that if $\msk_\vr
\setminus\{0\}\not=\emptyset$ then $c_\vr
=\inf\big\{\widetilde{\Phi}_\vr (u):\ u\in\msk_\vr
\setminus\{0\}\big\}$ (see an argument of \cite{Ding2008}). Using the
same iterative argument of \cite{Sere2} one obtains easily the
following

\begin{Lem}\label{critical points in W1,q}
Consider $\lm>0$ being a constant,
if $u\in\msk_\vr$ with $|\widetilde{\Phi}_\vr (u)|\leq C$,
then, for any $q\in[2,+\infty)$, $u\in
W^{1,q}(\R^3,\C^4)$ with
$\norm{u}_{W^{1,q}}\leq\Lam_q$ where $\Lam_q$ depends only on
$C$ and $q$.
\end{Lem}

\begin{proof}
See \cite{Sere2}. We outline the proof as follows. Firstly, from (\ref{D22}),
we write
\[
\begin{aligned}
u&=H_\om^{-1}\Big( \lm V_u \bt u
+ W_\vr (x)g(|u|)u  \Big).
\end{aligned}
\]
Now let $\rho:[0,\infty)\to[0,1]$
be a smooth function satisfying $\rho(s)=1$ if $s\in[0,1]$ and
$\rho(s)=0$ if $s\in[2,\infty)$. Then we have
\[
\begin{aligned}
g(s):=&\,g_1(s)+g_2(s)\\
=&\,\rho(s)g(s)+(1-\rho(s))g(s).
\end{aligned}
\]
Consequently,
$u=u_1+u_2+u_3$ with
\[
\begin{aligned}
u_1=&H_\om^{-1}\big( W_\vr \cdot g_1(|u|)u \big),\\
u_2=& \lm H_\om^{-1}\big( V_u \bt u \big),\\
u_3=&H_0^{-1}\big( W_\vr \cdot g_2(\jdz{u})u \big).
\end{aligned}
\]

Next we remark that, by H\"{o}lder's inequality, for $q\geq2$
$$
\big| V_u \bt u \big|_s\leq \big| V_u \big|_6 \cdot |u|_q
$$
with $\frac{1}{s}=\frac{1}{6}+\frac{1}{q}$ and, jointly with
\eqref{g-sigma0 estimate},
$$
\big|W_\vr
\cdot g_2(|u|)u\big|_t\leq C_1 |W|_\infty |u|_{t(p-1)}^{p-1},
$$
where $C_1>0$ is a constant.
Hence, we obtain
$$
u_1\in W^{1,2}\cap W^{1,3},\ u_2\in W^{1,s},\ u_3\in W^{1,t}.
$$
Then, denoting $s^*=\frac{3s}{3-s}$ and $t^*=\frac{3t}{3-t}$,
one sees $u\in W^{1,q}$ with $q=\min\{s^*,t^*\}$.

Starting with $q=2$, a standard bootstrap argument shows
that $u\in\cap_{q\geq2}L^q$,
$u_1\in\cap_{q\geq2}W^{1,q}$, $u_2\in\cap_{6>q\geq2}W^{1,q}$ and
$u_3\in\cap_{q\geq2}W^{1,q}$.

By Sobolev embedding theorems, $u\in C^{0,\ga}$ for some
$\ga\in(0,1)$. This, together with elliptic regularity (see
\cite{Trudinger}), shows $V_u\in W^{2,2}_{loc}(\R^3)\cap L^2(\R^3)$ and
$$
\|V_u\|_{W^{2,2}(B_1(x))}\leq C_2 \Big(\lm
|u|_{L^{4}(B_2(x))}^2 +\|V_u\|_{H^1(B_2(x))}\Big)
$$
for all $x\in\R^3$, with $C_2$ independent of $x$ and
$\vr $, where $B_r(x)=\{y\in\mathbb{R}^3:\jdz{y-x}<r\}$ for
$r>0$. Since $W^{2,2}(B_1(x))\hookrightarrow C^{0,\de}({B_1(x)})$,
$\de\in(0,\frac12)$, we have
\begin{\equ}\label{Ak local estimates}
\|V_u\|_{C\,^{0,\de}(B_1(x))}\leq C_3\Big(\lm
|u|_{L^{4}(B_2(x))}^2 +\|V_u\|_{H^1(B_2(x))}\Big)
\end{\equ}
for all $x\in\mathbb{R}^3$ with $C_3$ independent of $x$ and
$\vr $. Consequently $V_u \in L^\infty$, and this yields
$$
\jdz{V_u \bt u}_s\leq
\jdz{V_u}_\infty\jdz{u}_s.
$$
Thus $u_2\in\cap_{q\geq2}W^{1,q}$, and combining with
$u_1,u_3\in\cap_{q\geq2}W^{1,q}$ the conclusion is obtained.
\end{proof}

\begin{Rem}
Let $\msl_\vr$ denote the set of all least energy solutions
of $\widetilde{\Phi}_\vr $. If $u\in\msl_\vr$,
then $\widetilde{\Phi}_\vr (u)=c_\vr \leq C_{e_0}$. Recall that
$\msl_\vr$ is bounded in $E$ with upper bound $\Lam$
independent of $\vr $. Therefore, as a consequence of Lemma
\ref{critical points in W1,q} we see that, for each
$q\in[2,+\infty)$ there is $C_q>0$ independent of $\vr $ such that
\begin{\equ}\label{Cq estimate of L}
\|u\|_{W^{1,q}}\leq C_q\ \ \ \ \mathrm{for\ all}\
u\in\msl_\vr.
\end{\equ}
This, together with the Sobolev embedding theorem, implies that
there is $C_\infty>0$ independent of
$\vr $ with
\begin{\equ}\label{C-infty estimate of L}
\|u\|_{\infty}\leq C_\infty\ \ \ \ \mathrm{for\ all}\
u\in\msl_\vr.
\end{\equ}
\end{Rem}

\section{Proof of the main result}\label{PMT}

Throughout this section we assume $\om\in(-a,a)$, $(P_0)$ and
$(G_1)$-$(G_2)$ are satisfied. We also suppose, without loss
of generality, that $0\in\msc$.
The proof of the main theorem will be achieved in
three parts: {\it Existence, Concentration, and Exponential decay}.

%In order to give a better look at our work, in the sequel, we will
%state the Lemmas before the proofs. This enable us to describe the
%outline of our abstract theorem.

\subsection*{Part 1. Existence}
Keeping the notation of Section \ref{Preliminary results} we now turn
to the existence result of the main theorem.
Its proof is carried out in three lemmas. The modified problem gives
us an access to Lemma \ref{d to gamma}, which is the key ingredient
for Lemma \ref{d-eps is attained}.

Recall that $\tilde{\ga}_m$ denotes the least energy of
$\widetilde{\mst}_m$ (see the
subsection \ref{modified problem}), where
$\mu=m:=\max_{x\in\R^3}W(x)$, and $J_m$ denotes the associated reduction
functional on $E^+$. We remark that, since $0\in\msc$,
$W_\vr(x)\to W(0)=m$ uniformly on bounded sets of $x$. Our existence
results present as follows:

\begin{Lem}\label{d to gamma}
$c_\vr \to\tilde{\ga}_m$ as $\vr \to0$.
\end{Lem}

\begin{Lem}\label{d-eps is attained}
$c_\vr $ is attained for all small $\vr >0$.
\end{Lem}

\begin{Lem}\label{least energy solution compact}
$\msl_\vr $ is compact in $H^1(\R^3,\C^4)$, for
all small $\vr >0$.
\end{Lem}

\begin{proof}[Proof of Lemma \ref{d to gamma}] Firstly we show that
\begin{\equ}\label{d-eps>gamma}
\liminf\limits_{\vr \to0}c_\vr \geq\tilde{\ga}_m.
\end{\equ}
Arguing indirectly, assume that $\liminf_{\vr \to0}\,c_\vr
<\tilde{\ga}_m$. By the definition of $c_\vr$ and Proposition \ref{d-eps=c-eps}
we can choose an $e_j\in \msn_\vr$ and $\delta>0$ such that
$$
\max_{u\in E_{e_j}}\widetilde{\Phi}_{\vr _j}(u)\leq\tilde{\ga}_m-\delta
$$
as $\vr _j\to0$. Since $W_\vr (x)\leq m$, the representations of
$\widetilde{\Phi}_{\vr}$ and $\widetilde{\mst}_m$ imply that $\widetilde{\Phi}_\vr
(u)\geq \mst_m(u)$ for all $u\in E$ and $\vr$ small. Note also
that $\tilde{\ga}_m\leq J_m(e_j) \leq \max_{u\in E_{e_j}}\widetilde{\mst}_m(u)$.
Therefore we get, for all $\vr_j$ small,
$$
\tilde{\ga}_m-\delta\geq\max_{u\in E_{e_j}} \widetilde{\Phi}_{\vr _j}(u)
\geq\max_{u\in E_{e_j}}\widetilde{\mst}_m(u)\geq\tilde{\ga}_m,
$$
a contradiction.

\medskip

We now turn to prove the desired conclusion. Set $W^0(x)=m-W(x)$ and
$W^0_\vr (x)=W^0(\vr  x)$. Then
\begin{\equ}\label{tilde-Phi-eps----Gm}
\widetilde{\Phi}_\vr (u)=\widetilde{\mst}_m(u)+\int W^0_\vr
(x)G(|u|).
\end{\equ}

In virtue of Lemma \ref{erfle}, let $u=u^++u^-\in \msk_m$
such that $\widetilde{\mst}_m(u)=\tilde{\ga}_m$ and set
$e=u^+$. Surely, $e\in\msm_m$, $\msj_m(e)=u^-$ and
$J_m(e)=\tilde{\ga}_m$. There is a unique $t_\vr >0$ such that $t_\vr
e\in\msn_\vr $ and one has
\begin{\equ}\label{d-eps<I-eps(t-eps)}
c_\vr \leq I_\vr (t_\vr  e).
\end{\equ}
By Lemma \ref{For any e  Te} $t_\vr $ is bounded. Hence, without
loss of generality we can assume $t_\vr \to t_0$ as $\vr \to0$.
Using (\ref{the equa 1}), we infer
\[
\begin{aligned}
&\frac{1}{2}\|v_\vr\|^2+(I)
=\, \widetilde{\Phi}_\vr (w_\vr)-\widetilde{\Phi}_\vr (z_\vr)\\
=&\, \widetilde{\mst}_m(w_\vr)-\widetilde{\mst}_m(z_\vr)
 +\int W^0_\vr (x)\big( G(\jdz{w_\vr}) - G(\jdz{z_\vr}) \big)
\end{aligned}
\]
where, setting
\[
z_\vr=t_\vr e+\msj_m(t_\vr e), \  w_\vr=t_\vr e+h_\vr
(t_\vr e), \  v_\vr =z_\vr-w_\vr,
\]
\[
(I):=\int_0^1(1-s)
\big( \msf_\lm ''(w_\vr+sv_\vr )[v_\vr ,v_\vr ]
+\Psi_\vr ''(w_\vr+sv_\vr )[v_\vr ,v_\vr
]\big) ds.
\]
Taking into account that
\[
\begin{aligned}
&\int W^0_\vr (x)\big(G(|w_\vr|)-G(|z_\vr|)\big)\\
=&-\Re \int W^0_\vr (x)g(|z_\vr|)z_\vr \cdot\overline{v_\vr }
+\int_0^1(1-s)\msg_m''(z_\vr-sv_\vr )[v_\vr ,v_\vr ] \, ds\\
&-\int_0^1(1-s)\Psi_\vr ''(z_\vr-sv_\vr )[v_\vr
,v_\vr ] \, ds   \, ,
\end{aligned}
\]
setting
\[
\begin{aligned}
(II):=&\int_0^1(1-s)\Psi_\vr ''(z_\vr-sv_\vr
)[v_\vr ,v_\vr ] \,ds  \, ,
\end{aligned}
\]
following Remark \ref{auto inequ X}, one has
\[
\frac{1}{2}\norm{v_\vr }^2+(I)+(II)\leq
-\Re\int W^0_\vr (x)g(|z_\vr|) z_\vr \cdot\overline{v_\vr } \ .
\]
By noticing that $0\leq P^0_\vr (x)\leq m$, $(II)\geq0$ and
\[
\big| \msf_\lm ''(w_\vr+sv_\vr)[v_\vr ,v_\vr ] \big|
\leq \frac12 \|v_\vr\|^2  \, ,
\]
we deduce that
\begin{\equ}\label{v-eps estimate}
\frac{1}{4}\|v_\vr \|^2 \leq \int W^0_\vr (x)g(|z_\vr|)
 |z_\vr| \cdot |v_\vr| \, .
\end{\equ}

Since $t_\vr\to t_0$, it is clear that $\{z_\vr\}$, $\{w_\vr\}$
and $\{v_\vr\}$ are bounded and, particularly, for $q\in[2,3]$
$$
\limsup_{r\to\infty}\int_{|x|>r}|z_\vr|^q=0.
$$
Now we infer
\[
\begin{aligned}
&\int \bkt{W^0_\vr (x)}^{q/(q-1)} |u_\vr|^q  \\
=&\bigg(\int_{|x|\leq r}+\int_{|x|>r}\bigg){W^0_\vr (x)}^{q/(q-1)} |u_\vr|^q \\
\leq&\int_{|x|\leq r}\bkt{W^0_\vr (x)}^{q/(q-1)} |u_\vr|^q
+m^{q/(q-1)}\int_{|x|>r} |u_\vr|^q  \\
=&\ o(1)
\end{aligned}
\]
as $\vr \to0$. Thus by (\ref{v-eps estimate}) one has
$\norm{v_\vr}^2\to0$, that is, $h_\vr
(t_\vr  e)\to\msj_m(t_0 e)$. Consequently,
$$
\int W^0_\vr (x)G(\jdz{w_\vr})\to0
$$
as $\vr \to0$. This, jointly with (\ref{tilde-Phi-eps----Gm}), shows
\[
\widetilde{\Phi}_\vr (w_\vr)=\widetilde{\mst}_m(w_\vr)+o(1)
 =\widetilde{\mst}_m(z_\vr)+o(1),
\]
that is,
$$
I_\vr (t_\vr  e)=J_m(t_0e)+o(1)
$$
as $\vr \to0$. Then, since
$$
J_m(t_0e)\leq\max_{v\in E_e}\widetilde{\mst}_m(v)=J_m(e)=\tilde{\ga}_m,
$$
we obtain by using \eqref{d-eps>gamma} and
(\ref{d-eps<I-eps(t-eps)})
$$
\tilde{\ga}_m\leq
\lim_{\vr \to0}c_\vr
\leq\lim_{\vr \to0}I_\vr (t_\vr
e)=J_m(t_0e)\leq\tilde{\ga}_m,
$$
hence, $c_\vr \to\tilde{\ga}_m$.
\end{proof}

\medskip

\begin{proof}[Proof of Lemma \ref{d-eps is attained}]
Given $\vr >0$, let $\{u_n\}\subset\msn_\vr $ be a
minimization sequence: $I_\vr (u_n)\to c_\vr $. By the Ekeland
variational principle we can assume that $\{u_n\}$ is in fact a
$(PS)_{c_\vr }$-sequence for $I_\vr $ on $E^+$ (see
\cite{Pankov,Willem}). Then $w_n=u_n+h_\vr (u_n)$ is a $(PS)_{c_\vr
}$-sequence for $\widetilde{\Phi}_\vr $ on $E$. It is clear that
$\{w_n\}$ is bounded, hence is a $(C)_{c_\vr}$-sequence. We can assume
without loss of generality that $w_n\rightharpoonup w_\vr =w_\vr
^++w_\vr ^-\in\msk_\vr $ in $E$. If $w_\vr \not=0$ then
$\widetilde{\Phi}_\vr (w_\vr )=c_\vr $. So we are going to show that
$w_\vr \not=0$ for all small $\vr >0$.

To this end, take $\limsup_{|x|\to\infty}W(x)<\ka<m$ and
define
$$
W^\kappa(x)=\min\{\kappa,W(x)\}.
$$
Set $A:=\{x\in\R^3:W (x)>\ka\}$ and $A_\vr:=\{x\in\R^3:\vr x\in A\}$.
Following $(P_0)$, $A_\vr$ is a bounded set for any fixed $\vr$.
Consider the functional
$$
\widetilde{\Phi}_\vr
^\ka(u)=\frac{1}{2}\big(\|u^+\|^2-\|u^-\|^2 \big)
-\msf_\lm (u)-\int W^\ka_\vr (x)G(|u|)
$$
and as before define correspondingly $h_\vr ^\ka:E^+\to E^-$,
$I_\vr ^\ka:E^+\to\R$, $\msn_\vr ^\ka$, $c_\vr
^\kappa$ and so on. As done in the proof of Lemma \ref{d to gamma},
\begin{\equ}\label{d-eps-sigma  to  gamma-sigma}
\lim_{\vr \to0}c_\vr ^\ka=\tilde{\ga}_\ka.
\end{\equ}

Assume by contradiction that there is a sequence $\vr _j\to0$ with
$w_{\vr _j}=0$. Then $w_n=u_n+h_{\vr _j}(u_n)\rightharpoonup0$ in
$E$, $u_n\to0$ in $L_{loc}^q$ for $q\in[1,3)$, and $w_n(x)\to0$ a.e.
in $x\in\R^3$. Let $t_n>0$ be such that
$t_nu_n\in\msn_{\vr _j}^\ka$. Since $u_n\in\msn_\vr
$, it is not difficult to see $\{t_n\}$ is bounded and one may
assume $t_n\to t_0$ as $n\to\infty$.  Remark that
$h_{\vr _j}^\ka(t_nu_n)\rightharpoonup0$ in $E$ and $h_{\vr
_j}^\ka(t_nu_n)\to0$ in $L_{loc}^q$ for $q\in[1,3)$ as
$n\to\infty$ (see \cite{Ackermann}). Moreover, we remind that
$$
\widetilde{\Phi}_{\vr_j}(t_n u_n + h_{\vr_j}^\ka(t_n u_n))\leq
I_{\vr _j}(t_n u_n) \leq
I_{\vr _j}(u_n).
$$
So, we obtain
\[
\begin{aligned}
c_{\vr _j}^\kappa&\leq I_{\vr _j}^\kappa(t_nu_n)
=\widetilde{\Phi}_{\vr _j}^\kappa
(t_nu_n+h_{\vr _j}^\kappa(t_nu_n))\\
 &=\widetilde{\Phi}_{\vr _j}(t_nu_n+h_{\vr _j}^\kappa(t_nu_n))
 +\int \big(P_{\vr _j}(x)-P_{\vr _j}^\kappa(x)\big)
 G\big(|t_nu_n+h_{\vr _j}^\kappa(t_nu_n)|\big)\\
 &\leq I_{\vr _j}(u_n)+\int_{A_{\vr _j}}
 \big(P_{\vr _j}(x)-P_{\vr _j}^\kappa(x)\big)
 G\big(|t_nu_n+h_{\vr _j}^\kappa(t_nu_n)|\big)\\
 &=c_{\vr _j}+o(1)
\end{aligned}
\]
as $n\to\infty$. Hence $c_{\vr _j}^\kappa\leq c_{\vr _j}$. By
(\ref{d-eps-sigma  to  gamma-sigma}), letting $j\to\infty$ yields
$$
\tilde{\ga}_\ka \leq \tilde{\ga}_m,
$$
which contradicts $\tilde{\ga}_m < \tilde{\ga}_\ka$.
\end{proof}

\medskip

\begin{proof}[Proof of Lemma \ref{least energy solution compact}]
Since $\msl_\vr \subset B_\Lam$ for all small $\vr >0$,
assume by contradiction that, for some $\vr _j\to0$,
$\msl_{\vr _j}$ is not compact in $E$. Let
$u_n^j\in\msl_{\vr _j}$ with $u_n^j\rightharpoonup0$ as
$n\to\infty$. As done in proving the Lemma \ref{d-eps is attained},
one gets a contradiction.

Now let $\{u_n\}\subset\msl_\vr $ such that $u_n\to u$ in
$E$. We recall that $H_\om=i\al\cdot\nabla-a\bt + \om$, by
$$
H_\om u_n= \lm V_{u_n} \bt u_n
+ W_\vr (x)g(|u_n|)u_n
$$
and
$$
H_\om u= \lm V_u \bt u
+ W_\vr (x)g(|u|)u
$$
we deduce
\begin{\equ}\label{H1 estimate}
\begin{split}
\jdz{H_\om(u_n-u)}_2\leq \lm
 \big| V_{u_n} u_n-V_u u \big|_2
 +\big| W_\vr \cdot \big( g(\jdz{u_n})u_n-g(\jdz{u})u \big)
 \big|_2  \, .
\end{split}
\end{\equ}
%And a standard calculus shows that
%\[
%\begin{aligned}
%\big|V_{u_n} u_n -V_u u \big|_2 \leq &\,
%|u_n|_\infty^{1/6}
%\big| V_{u_n}-V_u \big|_6 |u_n|_{5/2}^{5/6}\\
% &+ |u_n-u|_\infty^{1/6}\big| V_u \big|_6 |u_n-u|_{5/2}^{5/6}
%\end{aligned}
%\]
%and
%\[
%\begin{aligned}
%&\big| W_\vr \cdot \big( g(|u_n|)u_n-g(|u|)u \big) \big|_2\\
%\leq& |W|_\infty \cdot \big| g(|u_n|)-g(|u|) \big|_\infty^{\frac12}
%\cdot
%\big| \big( g(|u_n|)-g(|u|) \big)^{1/2}u_n \big|_2\\
% &+|W|_\infty \cdot \big| g(\jdz{u}) \big|_\infty
% \cdot |u_n-u|_2.
%\end{aligned}
%\]
Invoking Lemma \ref{lemma of Aj} and $u_n\to u$ in
$L^q(\R^3,\C^4)$ for all $q\in[2,3]$, one gets
$\jdz{H_\om (u_n-u)}_2\to0$ as $n\to\infty$, and that is, $u_n\to u$ in
$H^1(\R^3,\C^4)$.
\end{proof}

\subsection*{Part 2. Concentration}
It is contained in the following lemma. To prove the lemma, it
suffices to show that for any sequence $\vr_j\to0$ the corresponding
sequence of solutions $u_j\in\msl_{\vr_j}$ converges, up to a
shift of $x$-variable, to a least energy solution of the limit
problem \eqref{the limit problem}.

\begin{Lem}\label{concentration}
Suppose that $\nabla W$ is bounded.
There is a maximum point $x_\vr $ of $|u_\vr|$ such that
$\dist(y_\vr ,\msc)\to0$ where $y_\vr =\vr  x_\vr $,
and for any such $x_\vr $, $v_\vr (x):=u_\vr (x+x_\vr )$ converges
to a ground state solution of (\ref{the limit problem}) in $H^1$
as $\vr \to0$.
\end{Lem}

\begin{proof}
Let $\vr _j\to0$, $u_j\in\msl_j$, where
$\msl_j=\msl_{\vr _j}$. Then $\{u_j\}$ is bounded. A
standard concentration argument (see \cite{Lions}) shows that there
exist a sequence $\{x_j\}\subset\R^3$ and constant $R>0$,
$\de>0$ such that
$$\liminf_{j\to\infty}\int_{B(x_j,R)}|u_j|^2 \geq \de.$$
Set
$$
v_j=u_j(x+x_j),
$$
and denoted by
$\hat{W}_j(x)=W(\vr _j(x+x_j))$,
one easily checks that $v_j$ solves
\begin{\equ}\label{vj equa}
H_\om v_j - \lm V_{v_j}\bt v_j
=\hat{W}_j \cdot g(|v_j|)v_j,
\end{\equ}
with energy
\[
\begin{aligned}
S(v_j)&:=\frac{1}{2}\big( \|v_j^+\|^2-\|v_j^-\|^2 \big)
- \Ga _\lm (v_j)-\int\hat{W}_j(x)G(|v_j|)  \\
 &=\widetilde{\Phi}_j(v_j)=\Phi_j(v_j)= \Ga _\lm(v_j)
 +\int\hat{W}_j(x)\widehat{G}(|v_j|)  \\
 &=c_{\vr _j}.
\end{aligned}
\]
Additionally, $v_j\rightharpoonup v$ in $E$ and $v_j\to v$ in
$L_{loc}^q$ for $q\in[1,3)$.

We now turn to prove that $\{\vr _jx_j\}$ is bounded. Arguing
indirectly we assume $\vr _j\jdz{x_j}\to\infty$ and get a
contradiction.

Without loss of generality assume $W(\vr _jx_j)\to W_\infty$.
By the boundness of $\nabla W$, one sees that $\hat{W}_j(x)\to W_\infty$
uniformly on bounded sets of $x$.
Surely, $m>W_\infty$ by $(P_0)$. Since for any $\psi\in C_c^\infty$
\[
\begin{aligned}
0&=\lim_{j\to\infty}\int\big(H_\om v_j
-\lm  V _{v_j} \bt v_j
-\hat{W}_j g(|v_j|)v_j \big) \bar \psi \\
 &=\lim_{j\to\infty}
 \int \big( H_\om v - \lm V_v \bt v - W_\infty g(|v|)v
 \big) \bar \psi ,
\end{aligned}
\]
hence $v$ solves
$$
i\alpha\cdot\nabla v-a\bt v + \om v - \lm V_v\bt v = W_\infty g(|v|)v.
$$
Therefore,
$$
S_\infty(v):=\frac{1}{2}\big( \|v^+\|^2-\|v^-\|^2 \big)
-\Ga_\lm(v)
-{W}_\infty \int G(|v|) \geq \tilde{\ga}_{W_\infty}.
$$
It follows from $m>P_\infty$, by
Proposition \ref{d-eps=c-eps}, one has
$\tilde{\ga}_m < \tilde{\ga}_{W_\infty}$.
Moreover, by the Fatou's lemma,
\[
\begin{aligned}
&\lim_{j\to\infty}\int\hat{W}_j(x)\widehat G (|v_j|)
\geq \int{W}_\infty\widehat G (|v|) \, .
\end{aligned}
\]
Consequently, noting that $\liminf_{j\to\infty}\Ga_\lm(v_j)\geq \Ga_\lm(v)$,
we have
$$
\tilde{\ga}_m<\tilde{\ga}_{W_\infty} \leq
S_\infty(v) \leq \lim_{j\to\infty}c_{\vr_j}=\tilde{\ga}_m  \, ,
$$
a contradiction.

Thus $\{\vr _jx_j\}$ is bounded. And hence, we can assume
$y_j=\vr _jx_j\to y_0$. Then $v$ solves
\begin{\equ}\label{v solves equa}
i\alpha\cdot\nabla v-a\beta v + \om v - \lm V_v\bt v=W(y_0)g(|v|)v \, .
\end{\equ}
Since $W(y_0)\leq m$, we obtain
$$
S_0(v):=\frac{1}{2}\big( \|v^+\|^2-\|v^-\|^2
\big) -\Ga_\lm(v) - W(y_0)\int G(|v|)
\geq \tilde{\ga} _{W(y_0)}\geq \tilde{\ga}_m.
$$
Again, by Fatou's lemma, we have
$$
S_0(v)=\int{P}(y_0)\widehat G (\jdz{v}) + \Ga_\lm(v)
\leq\lim_{j\to\infty}c_{\vr _j}=\tilde{\ga}_m.
$$
Therefore, $\ga_{P(y_0)}=\ga_m$, which implies
$y_0\in\msc$ by Proposition \ref{d-eps=c-eps}. By virtue of Lemma
\ref{critical points in W1,q} and (\ref{C-infty estimate of L}) it
is clear that one may assume that $x_j\in\mathbb{R}^3$ is a maximum
point of $\jdz{u_j}$. Moreover, from the above argument we readily
see that, any sequence of such points satisfies
$y_j=\vr _jx_j$ converging to some point in $\msc$ as
$j\to\infty$.

In order to prove $v_j\to v$ in $E$, recall that as the argument shows
$$\lim_{j\to\infty}\int\hat{W}_j(x)\widehat G (\jdz{v_j})
=\int{W}(y_0)\widehat G (\jdz{v}).$$
By $(G_2)$ and the decay of $v$, using the Brezis-Lieb lemma,
one obtains $|v_j-v|_\sa\to0$, then
$|v_j^\pm-v^\pm|_\sa\to0$ by (\ref{lpdec}). Denote $z_j=v_j-v$.
Remark that $\{z_j\}$ is bounded in $E$ and $z_j\to0$ in
$L^\sigma$, therefore $z_j\to0$ in $L^q$ for all $q\in(2,3)$. The scalar
product of (\ref{vj equa}) with $z_j^+$ yields
$$
\big\langle v_j^+, z_j^+\big\rangle=o(1).
$$
Similarly, using the decay of $v$ together with the fact
that $z_j^\pm\to0$ in $L_{loc}^q$ for $q\in[1,3)$, it follows from
(\ref{v solves equa}) that
$$
\big\langle v^+, z_j^+\big\rangle=o(1).
$$
Thus
$$
\|z_j^+\|=o(1),
$$
and the same arguments show
$$
\|z_j^-\|=o(1),
$$
we then get $v_j\to v$ in $E$, and the arguments in Lemma
\ref{least energy solution compact} show that $v_j\to v$ in
$H^1$.
\end{proof}

\subsection*{Part 3. Exponential decay}
See the following Proposition \ref{exp decay}.
%The first is to show a
%priori decay for the sequence of the least energy solutions obtained
%in Lemma \ref{d-eps is attained}. This we shall do with the help of
%the sub-solution estimate to show that $\{v_\vr\}$ is uniformly
%small at the infinity.
For the later use denote $D=i\al\cdot\nabla$ and, for
$u\in\msl_\vr $,  rewrite \eqref{D22} as
$$
Du=a\bt u - \om u+ \lm V_u \bt u
 + W_\vr (x)g(|u|)u.
$$
Acting the operator $D$ on the two sides and noting that
$D^2=-\De$, we get
\begin{\equ}\label{elliptic equ}
\aligned
\De u=&\, \big(a+\lm V_u \big)^2u -
\big( \om  - W_\vr \cdot g(|u|) \big)^2u\\
 &- D\big( \lm V_u + W_\vr \cdot g(|u|)\big) u \, .
\endaligned
\end{\equ}

Now define
\[
\sgn\, u = \left\{
\aligned
& \frac{\bar u}{|u|} \quad {\rm if\ } u\neq0 \, , \\
& 0 \quad {\rm if\ } u=0 \, .
\endaligned \right.
\]
By Kato's inequality \cite{R. Dautray}, there holds
\[
\De |u| \geq \Re \big[ \De u (\sgn\, u)  \big].
\]
Note that
\[
\Re \big[
D\big( \lm V_u + W_\vr \cdot g(|u|)\big) u (\sgn\,u)
\big] = 0  \, .
\]
Then, we obtain
\begin{\equ}\label{XXX}
\De |u|\geq \big(a+\lm V_u \big)^2 |u| -
\big( \omega  - W_\vr \cdot g(|u|) \big)^2 |u| \, .
\end{\equ}
To get the uniformly decay estimate for the semi-classical states,
we first need the following result:

\begin{Lem}\label{v-eps uniformly to 0}
Let $v_\vr $ and $V_{v_\vr }$ be given in the
proof of Lemma \ref{concentration}. Then $|v_\vr (x)|$ and
$|V_{v_\vr }(x)|$ vanish at infinity uniformly in
$\vr >0$ small.
\end{Lem}

Due to \eqref{XXX}, we remark that Lemma \ref{v-eps uniformly to 0} makes it feasible
to choose $R>0$ (independent of $\vr$) such that
\[
\De |v_\vr|\geq \frac{a^2-\om^2}{2} |v_\vr|
\quad {\rm for\ } |x|\geq R \, .
\]
And at this point, applying the maximum principle (see \cite{Rabier}),
we easily have

\begin{Prop}\label{exp decay}
Let $v_\vr\in E$ be given in the proof of Lemma \ref{concentration},
then $v_\vr$ exponentially decays at infinity uniformly in $\vr>0$ small.
More specifically, there exist $C,c>0$ independent of $\vr$ such that
\[
|v_\vr(x)|\leq C e^{-c|x|}  \, .
\]
Consequently, we infer that
\[
|u_\vr(x)|\leq C e^{-c|x-x_\vr|}  \, .
\]
\end{Prop}

Now, we turn to prove Lemma \ref{v-eps uniformly to 0}. To begin with,
we remind that \eqref{XXX} together with the regularity results for $u$ (see Lemma
\ref{critical points in W1,q})
implies there is $M>0$ (independent of $\vr$) satisfying
$$
\Delta |u| \geq -M |u|.
$$
It then follows from the sub-solution estimate \cite{Trudinger, Simon} that
\begin{\equ}\label{sub solu}
|u(x)|\leq  C_0 \int_{B_1(x)} |u(y)| dy
\end{\equ}
with $C_0>0$ independent of $x$, $\vr$ and $u\in\msl_\vr$.

\begin{proof}[Proof of Lemma \ref{v-eps uniformly to 0}]
Assume by contradiction that there exist $\de>0$ and
$x_\vr \in\R^3$ with $|x_\vr|\to\infty$ such that
$$
\de \leq |v_\vr(x_\vr)| \leq C_0 \int_{B_1(x_\vr)} |v_\vr(y)| dy \, .
$$
Since $v_\vr\to v$ in $E$, we
obtain, as $\vr\to0$,
\[
\begin{aligned}
\de &\leq C_0\Big(\int_{B_1(x_\vr)}|v_\vr|^2\Big)^{1/2}\\
 &\leq C_0\Big(\int |v_\vr-v|^2\Big)^{1/2}+C_0\Big(\int_{B_1(x_\vr)}|v|^2\Big)^{1/2}\to0,
\end{aligned}
\]
a contradiction. Now, jointly with (\ref{Ak local estimates}), one
sees also $|V_{v_\vr}(x)|\to0$ as $|x|\to\infty$
uniformly in $\vr >0$ small.
\end{proof}

With the above arguments, we are ready to prove Theorem
\ref{main theorem}.

\begin{proof}[Proof of Theorem \ref{main theorem}]
Going back to system (\ref{D2}), with the variable substitution:
$x\mapsto x/\vr $, Lemma \ref{d-eps is attained} jointly with Lemma
\ref{critical points in W1,q} and the elliptic regularity shows that,
for all $\vr >0$ small,
Eq.(\ref{D2}) has at least one ground state solution $(\va_\vr,\phi_\vr)\in
\cap_{q\geq2}W^{1,q}\times C^2$. Moreover,
by Lemma \ref{least energy solution compact} and Lemma \ref{lemma of Aj},
one easily checks the compactness of the ground states.

Assume additionally $\nabla W$ is bounded,
Lemma \ref{concentration} is nothing but the
concentration result.
And finally, Proposition \ref{exp decay} gives the decay estimate.
\end{proof}

\medskip

\noindent {\it Acknowledgment.} \ The authors wish to thank
the referee very much for his/her valuable comments and
suggestions.

The work was supported by the National Science Foundation of China
(NSFC11331010, 11171286) and the Beijing Center for Mathematics and
Information Interdisciplinary Sciences.

\end{document}